\documentclass[12pt]{article}
\usepackage[cp1251]{inputenc}
\usepackage[T2A]{fontenc}
\usepackage[english]{babel}

\usepackage{amsmath,amssymb,amsfonts,amsthm}

\begin{document}

\begin{center}

{\large {{\bf
Measures on a Hilbert space that are invariant with respect to shifts and orthogonal transformations
}}}
\end{center}

\begin{center}
Sakbaev V.Zh.\footnote{Steklov Mathematical Institute of Russian Academy of Sciences, fumi2003@mail.ru}
\end{center}


{\centerline {\bf Abstract}}

A finitely-additive measure $\lambda $ on an infinite-dimensional real Hilbert space $E$ which is 
invariant with respect to shifts and orthogonal transformations has been defined. This measure can be considered as the analog of the Lebesgue measure in the sense of its invariance with respect  to the above transformations. The constructed measure is defined on the ring $\cal R$ of subsets of the Hilbert space generated by measurable rectangles. A measurable rectangle is an infinite-dimensional parallelepiped such that the product of the lengths of its edges converges unconditionally.

This measure is obtained as a continuation of a family of shift-invariant measures $\lambda _{\cal E}$, where each measure $\lambda _{\cal E}$ is defined on the ring ${\cal R}_{\cal E}$ of measurable rectangles with edges collinear to the vectors of some orthonormal basis $\cal E$ in the space $E$. An equivalence relation is introduced on the set of orthonormal bases in terms of the transition matrix from one orthonormal basis to another. If two bases $\cal E$ and $\cal F$ are equivalent, then
the corresponding shift-invariant measures $\lambda _{\cal E}$ and $\lambda _{\cal F}$ coincide. If these two bases are not equivalent, then the restrictions of the measures $\lambda _{\cal E}$ and $\lambda _{\cal F}$ 
onto the intersection of their domains are equal to zero. This property allows to glue measures defined on the subset rings corresponding to different bases into the one measure $\lambda $ defined on the unique ring $\cal R$. The obtained measure $\lambda $ is invariant with respect to  shifts and rotations.
The decomposition of the measure $\lambda$ into the sum of mutually singular shift-invariant measures is obtained.

The paper describes the structure of the space $\cal H$ of numerical functions
square integrable with respect to the constructed shift and rotation-invariant measure $\lambda $. It is shown that the spaces ${\cal H}_{\cal E}$ of functions that are square integrable with respect to a measure $\lambda _{\cal E}$ associated with one equivalence class of bases form a closed subspace ${\cal H}_{\cal E}$ in the space $ \cal H$. As a result, there is the decomposition of the space $\cal H$ into the orthogonal sum of subspaces corresponding to all possible equivalence classes of bases. Unitary groups acting  by means of orthogonal transformations of the argument in the space $\cal H$ of square integrable functions are investigated.



\section{Introduction}

A nontrivial countably additive 
$\sigma $-finite  locally finite Borel left-invariant  measure on a topological group $G$ does not exist according to A. Weyl theorem if the group $G$ is not locally compact. 
Hence there is no  nontrivial countably additive 
$\sigma $-finite  locally finite Borel shift-invariant measure on an infinite-dimensional normed linear space.
Therefore, the studying of the shift-invariant measures on the Hilbert space deals with additive function of a set without some properties of Lebesgue measure. 
We study finitely-additive measures on an infinite-dimensional separable real Hilbert space so that these measures are
invariant under shifts to an arbitrary vector and arbitrary orthogonal transformations.
The focus of our research is  the spaces of functions that are square  integrable by an invariant measure. Unitary groups acting  by means of orthogonal transformations of the argument in the space of square integrable functions are investigated.

Thus, a shift-invariant measure on a topological group without the local compactness property (in particular a shift invariant measure on an infinite-dimensional linear space equipped with the norm or with the locally convex topology) is defined as following. A measure on a topological vector space is considered as an additive nonnegative function, which is defined on some ring of subsets of the space. However, this function of a set should have no at least one property of the Lebesgue measure listed in A. Weyl theorem 
\cite{Baker, Pantsulaya, s16, Vershik, Zavadski, SITO}. 

One of the approaches is based on the construction of countably additive measure without the $\sigma$-finiteness property.
The countably additive  measures on the topological vector spaces of numerical sequences  are introduced in \cite{Baker, Pantsulaya, Zavadski}. But the introduced measures are not $\sigma$-finite and locally finite.

In other approach, shift-invariant $\sigma$-finite locally finite measures on a separable Banach space are introduced in  \cite{s16, SITO}. But any of the constructed measures is not countably additive and is not defined on the ring of bounded Borel subsets.  
The paper \cite{SITO} describes the construction of a shift-invariant measure on a Hilbert space as an additive function of a set, which is defined on some ring of subsets of Hilbert space.
This ring of subsets (the domain of an additive function of a set) is not invariant with respect to any orthogonal transformation since it depends on the choice of orthonormal basis (ONB) in the Hilbert space.

The studying of the finitely-additive measure on a Hilbert space, which is invariant with respect to both shifts and rotation, is closed to the same problem in finite dimension euclidean space. The problem of the existence of an  invariant with respect to isometric transformation measure on a finite-dimensional euclidean space is investigated during last century in the following form.

{\it Does the measure $\lambda $ on the $d$-dimensional euclidean space exists so that this measure is

1) defined on any bounded subset of the euclidean space,

2) invariant with respect to any shift and any rotation,

3) normalized by the condition $\lambda ([0,1]^d)=1$?}

There is no countable additive measure with these properties for every $d\in \mathbb N$ according to the article by F. Hausdorff \cite{1914}.
In 1923 S. Banach had prove of existence of a  finitely-additive measure which is defined on a family of all bounded subsets of euclidean spaces ${\mathbb R}^d,\, d=1,2$,  so that this measure is invariant to any isometry
(see \cite{Natanson}, p. 81).  Hence  finitely-additive measures can admit invariance with respect to a more wide group than a countably additive one.

The paradox of Hausdorff-Banach-Tarskii gives some restriction onto the properties of a finitely-additive measure on the euclidean space which is discussed in  \cite{TTao}. In particular, there is no finitely-additive measure on the euclidean space ${\mathbb R}^d$ with dimension $d\geq 3$ which is defined on the ring of all subsets of this space and invariant with respect to shifts and rotations. However, according to work by S. Banach,  in the case $d=1,2$, there is a shift and rotation invariant finitely-additive nonnegative measure $\nu _d$ on the space ${\mathbb R}^d$ which is defined on every bounded subset of this space and is normalized by the condition $\nu _d([0,1)^d)=1$. 

In the class of countably additive measures, there is the only one normalized shift invariant complete Borel measure on the space ${\mathbb R}^d$ satisfying the normalization condition 3). This measure is the Lebesgue measure. It is additionally invariant with respect to rotation.

In 1923--1984 the following S. Ruziewicz problem
had been investigated.
Let ${\cal B}({\mathbb R}^n)$ be the ring of bounded Lebesgue measurable sets in the $n$-dimensional real space ${\mathbb R}^n$. Let $\lambda _n$ be the Lebesgue measure on ${\cal B}({\mathbb R}^n)$ normalized by $\lambda _n ([0,1)^n)=1$. {\it Is $\lambda _n$, up to proportionality, the unique finitely-additive isometry-invariant positive measure ${\cal B}({\mathbb R}^n)\to [0,+\infty )$? }

Banach gave a negative answer to this question for ${\mathbb R}^1$ and ${\mathbb R}^2$ in 1923. For ${\mathbb R}^n$ with $n\geq 3$ the question by Ruziewicz had been studied during a long time. The positive answer to this question for ${\mathbb R}^n$ with $n\geq 4$ was given independently and simultaneously in  \cite{Morg, Sull}. For ${\mathbb R}^3$ the positive answer to the Ruziewicz question was provided by G. Margulis \cite{Morg}. Also, in  \cite{ Morg, Sull, Drinfeld} it was proved that a finitely-additive measure on the unite sphere $S^{n-1}$ in the euclidean space ${\mathbb R}^n$ is proportional to Lebesgue measure on $S^{n-1}$ for every $n\geq 3$.

The Lebesgue measure can be defined as the complete shift-invariant extension of the measure which is defined on the ring $r({\mathbb R}^d)$ and is normalized by the condition 3). Here  $r({\mathbb R}^d)$ is the ring generated by the collection of bounded $d$-dimensional rectangles. Hence the  Ruziewicz problem can be reformulated in the following way: {\it what  finitely-additive isometry-invariant positive measure is the extension of  the measure which is defined on the ring $r({\mathbb R}^d)$ and is normalized by the condition 3)?}

However, in the separable infinite-dimensional real Hilbert space $E$ the Ruziewicz question should be reformulated  since there is no Lebesgue measure on the space $E$ or on the sphere in the space $E$.  In the present paper we study the following question. {\it What measure on the space $E$ exists  so that this measure is satisfied the conditions

1) it is invariant to any bijective isometric transformation of this space;
 
2) the domain of this measure is the ring $\cal R$ of subsets of the space $E$ which contains any measurable rectangles of the space $E$. A measurable rectangle in the space $E$ is an infinite-dimensional parallelepiped,
such that the product of the lengths of its edges converges unconditionally;

3) the normalized condition $\lambda (\{ x\in E:\ (x,e_k)\in [0,1)\ \forall \ k\in {\mathbb N}\})=1$ holds for some ONB ${\cal E}=\{ e_k\}$?}

Let us note that the ring $\cal R$ does not contain any bounded Borel set of the space $E$.

The measure with the mentioned invariance properties is important for studying the Hamiltonian flows in infinite-dimensional phase space. We obtain the integral invariant that is universal for the class of Hamiltonians.  In fact, the constructed measure is an invariant measure for class of  Hamiltonian flows in real Hilbert space with the standard symplectic structure, including a group of rotations and a group of shifts.

In the separable real euclidean space $\ell _2$, there is no normalized shift-invariant countably additive $\sigma$-finite and locally finite measure. There are different shift-invariant finitely-additive  $\sigma$-finite and locally finite measures, which is normalized by the condition $\lambda ([0,1)^{\mathbb N}\bigcap \ell _2)=1$  (see \cite{S-Smol, Zavadski}). In present paper, the extension of one of such measure onto the shift and rotation invariant normalized measure is introduced.

In this paper the measures on a Hilbert space combining the properties of shift-invariance and rotation-invariance have been introduced  by means of extension of finitely-additive shift-invariant measure on a Hilbert space
\cite{s16, SITO} onto the rotation invariant measure.

The existence of  a shift and rotation-invariant measure on a real separable Hilbert space is proved in \cite{s16}. This measure was introduced by using the axiom of choice. To achieve this aim a family $\{ \lambda _{\cal E}\}$ of shift-invariant measures is considered. Every measure $\lambda _{\cal E}$ depends on the choice of ONB $\cal E$ in the space $E$. The collection $\{ \lambda _{\cal E}\}$  of measures  is equipped with some total ordering. The shift and rotation-invariant measure $\lambda$ is defined as the extension of any of  shift-invariant measures $\lambda _{\cal E}$ by  using
of the transfinite induction procedure. 
The independence of the measure $\lambda $ of the choice of total ordering is proved in \cite {s16}.

It seems, that
another construction of a nontrivial rotation and shift-invariant measure on an infinite-dimensional Hilbert space has not be considered.
In particular, the extension of a countably additive shift-invariant measure on the space $\ell _2$ from \cite{Zavadski} onto the rotation and shift-invariant measure is still an unresolved problem.

Generalized shift invariant measure is investigated in  \cite{SSh} as the shift invariant functionals on the space of test functions of Schwartz type.
The constructed functionals have the properties of invariance with respect to the group of orthogonal mappings. But the question about the existence of a measure as the additive function of a set on some ring of subsets, connected with the constructed shift invariant functionals, is unresolved.

The description of deformation of a measure on the Hilbert space under the action of a linear mapping of this space is given by the infinite determinants in  \cite{NOMA}.
According to \cite{NOMA}, the linear mapping of Hilbert space with the finite determinant transforms the shift and rotation invariant measure $\lambda $ on this space into the measure that differs from the measure $\lambda$ by the multiplication on the absolute value of determinant of the linear mapping.

The paper \cite{Vershik} is devoted to the properties of a measure on a topological vector spaces which is invariant with respect on a vector from some admissible  subspace of the topological vector space. A constructed measure has every property of Lebesgue measures from the Weyl theorem except for the invariance with respect to a shift on arbitrary vector.

The invariance of a measure with respect to some group of unitary transformations is studied in   \cite{Pickerell}, \cite{Bufetov}. In these papers, the Pickrell measures on the space of infinite complex matrices and on the Grassman manifold of infinite-dimensional Hilbert space are constructed. Pickrell measures are the two-parametric family of probability measures on the space of complex matrices, which is invariant with respect to some infinite subgroup of  unitary operators group acting in the space of complex matrices by means of conjugation
\cite{BO}.

The measures on algebras of operators are studied  in \cite{BS19}. Some of these measures were defined by means of operator intervals  \cite{Bik18a}, \cite{BSu19}. The invariance of the introduced measures on algebras with respect to the action of some groups was obtained.

\bigskip

The purpose of this paper is to introduce a finitely-additive measure on an infinite-dimensional real Hilbert space so that this measure is invariant with respect to shift on any vector of a Hilbert space and with respect to any orthogonal mapping (i.e. the measure is shift  and rotation invariant). 
Moreover, the introduced measure is locally finite and $\sigma$-finite. But it is nether countably additive nor Borel measure.

Specifically,
the construction of a shift and rotation-invariant measure on a real separable Hilbert space $E$ is based on the analysis of the deformation under the action of orthogonal mappings on a  shift-invariant measure $\lambda _{\cal E}$ on the space $E$, where a measure $\lambda _{\cal E}$ depends on the choice of ONB $\cal E$  \cite{SITO}. 
We obtain the  criteria of an absolutely continuity of the measure  $\lambda _{\cal E}$  with respect to the image $\lambda _{\cal E}\circ {\bf U}$ of a measure $\lambda _{\cal E}$ under the action of the orthogonal mappings $\bf U$. The condition of the criteria of an absolute continuity of the measure $\lambda _{\cal E}\circ {\bf U}$ with respect to the measure $\lambda _{\cal E}$ is similar to  the condition of Poincare theorem (\cite{ArSad}, p. 399-400) for the existence of the  determinant of the matrix of the orthogonal mappings $\bf U$ in the ONB $\cal E$ (see the condition (\ref{crit}) below). If  the measure $\lambda _{\cal E}\circ {\bf U}$ is absolutely continuous with respect to the measure $\lambda _{\cal E}$,  then the measures $\lambda _{\cal E}\circ {\bf U}$ and $\lambda _{\cal E}$ coincide. In the opposite case,  the measures $\lambda _{\cal E}\circ {\bf U}$, $\lambda _{\cal E}$ are defined on the different rings ${\cal R}_{{\bf U}{\cal E}}$, ${\cal R}_{\cal E}$ respectively and the condition ${\cal R}_{{\bf U}{\cal E}}\bigcap {\cal R}_{\cal E}=\{ A\subset E:\ \lambda _{{\bf U}{\cal E}}(A)=\lambda _{\cal E}(A)=0\}$ holds. 

The equivalence relation $\sim $ on the set of ONB in the space $E$ is introduced in the following way. Two orthonormal bases $\cal E$ and $\cal F$ are equivalent to each other if and only if $\lambda _{\cal F}=\lambda _{\cal E}$. Let $\Sigma$ be a set of classes of equivalences $\{ {\cal E}\}$ of ONB.

This analysis of measures $\lambda _{\cal E}:\ {\cal R}_{\cal E}\to \infty $ with different ONB $\cal E$ gives the rule to define a shifts and rotations-invariant measure 
$\lambda :\ {\cal R}\to [0,+\infty )$, where $\cal R$ is the ring of subsets generated by the collection of sets $\bigcup\limits_{\cal E}{\cal R}_{\cal E}$.

The space ${\cal H}=L_2(E,{\cal R},\lambda , \mathbb C)$ of square integrable with respect shift and rotation invariant measure complex valued functions is studied.
The orthogonal decomposition 
$$
{\cal H}=\oplus _{\{{ \cal E}\}\in \Sigma} {\cal H}_{\{ {\cal E}\} }
$$ 
is obtained. Here ${\cal H}_{\{ {\cal E}\}}\equiv {\cal H}_{\cal E}$ for some ${\cal E}\in \{ {\cal E}\},\ \{ {\cal E}\}\in \Sigma$.
Any component ${\cal H}_{\{ {\cal E}\}}$ of orthogonal decomposition is invariant with respect to shift on any vector of the space $E$. The space $\cal H$ is invariant with respect to a shift and to an orthogonal transformation of the space $E$.

The structure of the present article is following.

Section 2 introduce the family shift-invariant measures $\{\lambda _{\cal E},\, {\cal E}\in {\cal S}\}$, where a measure $\lambda _{\cal E}:\ {\cal R}_{\cal E}\to [0,+\infty )$ depends on the choice of ONB $\cal E$ \cite{s16, BusS}.

Section 3 contains the description of mutual position for two ONB in the space $E$ in the terms of infinite orthogonal transition matrix. The condition of the proximity of one ONB to another one in the terms of transition matrix is introduced.
Section 4 exhibits  that if ONB $\cal E$ and $\cal F$ satisfy the proximity condition then the ring
${\cal R}_{\cal E}$ coincides with the ring ${\cal R}_{\cal F}$ and the measure $\lambda _{\cal E}$ coincides with the measure $\lambda _{\cal F}$. Section 5 demonstrates  that  if ONB $\cal E$ and $\cal F$ do not satisfy the proximity condition then $\lambda _{\cal E}(A)=0=\lambda _{\cal F}(A)\ \forall \ A\in {\cal R}_{\cal E}\bigcap {\cal R}_{\cal F}$.

To solve the problems of Section 4 and 5, the properties of the intersection of a measurable rectangle with its image under the action of a shift or an orthogonal mapping are studied.
This geometric problem is interesting as the infinite-dimensional generalization of the theory of
 $k$-dimensional sections of $n$-dimensional cubes  \cite{Ball, Ivanov}. The solving of this problem gives the possibility to introduces the equivalence relation on the set $\cal S$ of ONB of the space $E$ in the terms of the proximity condition from Section 3.

The proof of the existence of rotation invariant analog of Lebesgue measure on a Hilbert space is given in Section 6 without the using of the choice axiom. 
In Section 6 we prove that the orthogonal sum 
${\mathbb H}=\oplus _{\{ \cal E \}}{\cal H}_{\{ \cal E\} }$ on the set of equivalence classes of ONB is the Hilbert space of functions square integrable with respect to shift and rotation-invariant measure $\lambda :\ {\cal R}\to [0,+\infty )$.
The decomposition of the measure $\lambda$ into the sum of mutually singular shift-invariant measures is obtained.

In Section 7 we study the unitary group in the Hilbert space $\mathbb H$ which is generated by the orthogonal mapping of arguments of the functions from the space $\mathbb H$.
The condition of strong continuity  in the whole space $\mathbb H$ and the description of the continuity subspaces for this unitary groups are obtained in Section 7. These results are important for extending onto the infinite-dimensional case of the procedure of averaging of random orthogonal mappings and for obtaining the differential equation describing the mean values of compositions of independent random orthogonal mappings 
\cite{OSZ}.

In Section 8 the properties ob introduced invariant measure and its practical implication results are considered.

\section{Shift invariant measures on a Hilbert space
}

Let $E$ be a real separable Hilbert space.
Let $\cal S$ be a set of ONB in the space 
$E$.

Now we introduce a measure on the space $E$ which is invariant with respect to a shift on any vector of the space $E$ and which is connected with the choice of ONB
${\mathcal E}=\{ e_i\}$ (\cite{s16, SITO}).

A set $\Pi \subset E$ is called rectangle if there are an ONB ${\mathcal E}=\{ e_i\}$ and an elements $a,b\in \ell_{\infty}$ such that 
\begin{equation}\label{10}
\Pi =\{ x\in E:\ (x,e_j)\in [a_j,b_j)  \; \; \forall \, j\in {\mathbb{ N}}\}.
\end{equation}
A rectangle  (\ref{10}) is called measurable if either 
$\Pi =\varnothing $ or the following condition holds
\begin{equation}\label{11}
\sum\limits_{j=1}^{\infty}\max \{ 0, \ln (b_j-a_j)\} <\infty .
\end{equation}

Let $\mathcal K$ be a collection of measurable rectangles in the space $E$. Let $\rm r$ be a ring of subsets of the space $E$ generated by the collection 
$\mathcal K$. Let $\cal E$ be an ONB in the space $E$. Let ${\mathcal K}_{\mathcal E}$ be a set of measurable rectangles in the space $E$ such that the edges of any rectangle $\Pi \in {\mathcal K}_{\mathcal E}$ are collinear to the vectors of ONB   
$\mathcal E$. 
In other words, if $\Delta _j\subset \mathbb R$ is the projection of a set $\Pi$ onto the axis $Oe_j$ for any $j\in \mathbb N$ then $\Pi =\{ x\in E:\ (x,e_j)\in \Delta _j \ \forall \ j\in {\mathbb N}\}$.
Let ${\rm r}_{\mathcal E}$ be a ring of subsets of the space $E$ which is generated by the set ${\mathcal K}_{\mathcal E}$.

Let $\lambda $ be a function of a set such that the function $\lambda$ is defined on the collection of a sets $\mathcal K$ by the equality
\begin{equation}\label{12}
\lambda (\Pi )=\exp\Bigl[\sum\limits_{j=1}^{\infty }\ln (b_j-a_j)\Bigr] 
\end{equation}
for any nonempty measurable rectangle  (\ref{10}) and $\lambda (\varnothing)=0$. According to the condition (\ref{11}) $\lambda (\Pi )\in [0,+\infty )$ for any $\Pi \in \mathcal K$. Let  $\lambda _{\mathcal E}$ be the restriction of the function of a set $\lambda $ onto the collection of a sets ${\mathcal K}_{\mathcal E}$.

According to the works \cite{s16, SITO} the function 
$\lambda _{\mathcal E} $ is additive function on the collection of a sets ${\mathcal K }_{\mathcal E}$ and it has the unique extension onto the measure  $\lambda _{\mathcal E}:\ {\rm {r}} _{\mathcal E}\to \mathbb R$.  This measure $\lambda _{\mathcal E} :\ {\rm r} _{\mathcal E}\to [0,+\infty )$ is invariant with respect to a shift on any vector of the space  $E$.

A set $A\subset E$ is called $\lambda _{\mathcal E}$-measurable if for any $\epsilon >0$ there are sets $A_*,A^*\in {\rm r} _{\mathcal E}$ such that  $A_*\subset A\subset A^*$ and $\lambda (A^*\backslash A_*)<\epsilon $.
Then the collection ${\mathcal R} _{\mathcal E}$ of $\lambda _{\mathcal E}$-measurable subsets of the space $E$ is the ring. The measure   $\lambda _{\mathcal E}:\ {\rm r} _{\mathcal E}\to [0,+\infty )$ has the unique extension onto the ring ${\mathcal R} _{\mathcal E}$ by the equality $\lambda _{\mathcal E} (A)=\inf\limits_{A^*\in {\rm r} _{\mathcal E},\, A^*\supset A}\lambda (A^*)\;\; \forall \ A\in {\mathcal R} _{\mathcal E} $. 

The function of a set $\lambda _{\mathcal E} :{\mathcal R} _{\mathcal E} \to [0,+\infty )$ is the finitely-additive measure which is invariant with respect to shift on any vector of the space $E$ \cite{s16}. This measure is locally finite,  $\sigma$-finite, complete. But this measure is not $\sigma$-additive and it is not defined on the $\sigma$-ring of bounded Borel subsets. In particular, the ring ${\mathcal R} _{\mathcal E}$ does non contain a ball of the space $E$ with a sufficiently large radius \cite{Busovikov}.

Thus for a given ONB
${\cal E}\in \cal S$ there is the ring of subsets ${\cal R}_{\cal E}$ and there is the shift-invariant finitely-additive locally finite and $\sigma$-finite measure $\lambda _{\cal E}:\ {\cal R}_{\cal E}\to [0,+\infty )$.

The paper \cite{s16} describes the procedure of extension of the family of measures  $\lambda _{\cal F}:\ {\cal R}_{\cal F}\to [0,+\infty );\ {\cal F}\in {\cal S}$ onto the measure $\lambda :\ {\cal R}\to [0,+\infty )$ where ${\cal R}$ is the ring generated by the collection of sets $\bigcup\limits_{{\cal F}\in {\cal S}}{\cal R}_{\cal F}$ (therefore the ring $\cal R$ contains the collection of sets  $\cal K$).

The existence of the measure  $\lambda :\ {\cal R}\to [0,+\infty )$ such that $\lambda |_{{\cal R}_{\cal F}}=\lambda _{\cal F}\ \forall \ {\cal F}\in {\cal S}$ had been proved in the work \cite{s16} by using of some total ordering $\prec $ on the set $\cal S$ of ONB in the space $E$ and by applying of transfinite induction procedure.
It had been proved that the measure $\lambda$

1) is invariant with respect to any orthogonal mapping and to a shift on any vector ${\bf h}\in E$;

2) does not depend on the choice of total ordering $\prec$ on the set $\cal S$. 

The dependence of properties of the ring ${\cal R}_{\cal E}\bigcap {\cal R}_{\cal F}$ on the mutual position of two ONB 
$\cal E$ and $\cal F$ is not considered in the paper \cite{s16}. 

In present article we describe the dependence of the ring  ${\cal R}_{\cal E}\bigcap {\cal R}_{\cal F}$ on the mutual position of two ONB 
$\cal E$ and $\cal F$.

\section{Proximity for two ONB and the determinant of orthogonal transition matrix
}

Let us study the dependence of the ring ${\cal R}_{\cal EF}={\cal R}_{\cal E}\bigcap {\cal R}_{\cal E}$ on the  mutual position for two ONB ${\cal E},\,{\cal F}\in {\cal S}$ in the space $E$.
The description of this dependence gives the opportunity to define the procedure of the extension of the family of measures 
$\lambda _{\cal F},\ {\cal F}\in {\cal S},$ onto the rotation invariant measure $\lambda :\ {\cal R}\to [0,+\infty )$.

Let ${\bf U}$ be an orthogonal operator in the space $E$. Let ${\cal E},\, {\cal F}$ be a pair of ONB in the space $E$ such that ${\cal F}={\bf U}({\cal E})$. Let us consider the measures $\lambda _{\cal E}:\ {\cal R}_{\cal E}\to [0,+\infty )$ and $\lambda _{\cal F}:\ {\cal R}_{\cal F}\to [0,+\infty )$. We study the measures 
$\lambda _{\cal E}|_{{\cal R}_{\cal EF}}$ and $\lambda _{\cal 
F}|_{{\cal R}_{\cal EF}}$ where ${\cal R}_{\cal EF}={\cal R}_{\cal E}\bigcap {\cal R}_{\cal F}$.

Let $\|C\|=\|c_{i,j}\|$ be the matrix of transition of the basis $\cal E$ into the basis $\cal F$. Hence
the matrix elements  are given by the equalities
$c_{i,j}=(e_i,f_j)=(e_i,{\bf U}e_j),\ i,j\in {\mathbb N}$. Therefore $\sum\limits_{k=1}^{\infty }c_{k,i}c_{k,j}=\delta _{ij},\ i,j\in {\mathbb N}$ and $(e_k,{\bf U}^{-1}e_l)=c_{l,k},\ k,l\in \mathbb N$ where $\delta _{i,j}$ is the Kronecker symbol.

Since $\| c_{\cdot ,j}\|$ is the coordinates of the unit vector  
$f_j,\, j\in \mathbb N$, with respect to the basis  $\cal E$ then the sequence $\{ c_{\cdot ,j}\}$ is the unit vector of Hilbert space $\ell_2$. But the sequence  $\{ c_{\cdot ,j}\}$ can be not belong to the space $\ell_1$. We will show that if $\{ c_{\cdot ,j}\} \notin \ell _1$ then  $\lambda _{\cal E}(A)=0=\lambda _{\cal F}(A)$ for any set
$A\in {\cal R}_{\cal E}\bigcap {\cal R}_{\cal F}$. 
For a given ONB  $\cal E$ the symbol $L_1({\cal E})$ notes the linear subspace $L_1({\cal E})=\{ x\in E :\ \{ (x,e_j)\} \in \ell _1\}$.
If the conditions $f_j\in L_1({\cal E})\ \forall \ j\in \mathbb N$ and 
$e_j\in L_1({\cal F})\ \forall \ j\in \mathbb N$ hold then the property of absolute continuity of measures $\lambda _{\cal E},\ \lambda _{\cal F}$ with respect to each other is controlled by the following conditions on the pair of bases  
${\cal E},\, {\cal F}$
\begin{equation}
\label{crit}
\prod _{j=1}^{\infty}\| c_{\cdot ,j}\| _{\ell_1}<+\infty ;
\end{equation}
\begin{equation}
\label{crit'}
\prod _{j=1}^{\infty}\| c_{j,\cdot }\| _{\ell_1}<+\infty .
\end{equation}

We prove that the condition (\ref {crit}) is equivalent to the condition (\ref {crit'}). If the conditions   (\ref {crit}) and (\ref {crit'}) are satisfied then  ${\cal R}_{\cal E}\bigcap {\cal R}_{\cal F}={\cal R}_{\cal E}={\cal R}_{\cal F}$ and the equality  $\lambda _{\cal E}(A)=\lambda _{\cal F}(A)$ holds for any $A\in {\cal R}_{\cal E}\bigcap {\cal R}_{\cal F}$.  In the opposite case the measures  $\lambda _{\cal E}$ and $\lambda _{\cal F}$ take only zero values on an arbitrary set of the ring ${\cal R}_{\cal E}\bigcap {\cal R}_{\cal F}$.

The next lemma defines the geometrical sense of values of the products in left hand side of inequalities
(\ref{crit}) and (\ref{crit'}).

{\bf Lemma 3.1}. {\it Let $\Pi_{0,1}$ be a unit rectangle from the collection of sets ${\cal K}_{\cal E}$.
Then}
$$
\inf\limits_{Q\in {\cal K}_{\cal F}:\ \Pi _{0,1}\subset Q }\lambda _{\cal F}(Q)=\prod _{j=1}^{\infty}\| c_{\cdot ,j}\| _{\ell_1}.
$$

{\bf Proof}. According to the conditions $Q\in {\cal K}_{\cal F}:\ \Pi_{0,1}\subset Q$ the length $l_j$ of the orthogonal projection of the rectangle $Q$ on the line $Of_j=\{ x=tf_j,\, t\in {\mathbb R}\},\, j\in \mathbb N$ (i.e. the length of $j$-th edge of the rectangle $Q$) no less than the sum of the lengths of the orthogonal projections of edges of the rectangle $\Pi _{0,1}$ onto the line $Of_j$: 
 $l_j \geq \sum\limits_{i=1}^{\infty }|(e_i,f_j))|\ \forall \ j\in \mathbb N$. Therefore we obtain the statement of lemma 3.1.
\hfill$\Box$

{\bf Lemma 3.2}. {\it Let $\{ c_k\}\in \ell_2$ and $\| \{ c_k\} \|_{\ell_2}=1$.  If $\max\limits_{j\in {\bf N}}|c_j|=\alpha ,$  then $\|c _j\|_{\ell_1}\geq \alpha +{\sqrt {1-\alpha ^2}}$. In particular,  $\| \{ c_k\}\|_{\ell_1}\geq 1$.}

The statement is the consequence of the inequality  $|c_1|+\ldots+|c_m|\geq {\sqrt{ c_1^2+\ldots+c_m^2}}$ which holds for any $m\in \mathbb N$ and for any collection of complex numbers $c_1, \ldots, c_m$. \hfill$\Box$

{\bf Corollary 3.3.}  {\it Let $c^{(m)}$ be a sequence of vectors of the space $\ell_2$ with coordinates $c_k^{(m)},\, k\in \mathbb N$ such that $\lim\limits_{m\to \infty}\| \{ c_k^{(m)}\} \|_{{\ell}_2}=1$.  If $\lim\limits_{m\to \infty }\| c^{(m)}\|_{{\ell}_1}=1$ then  
$\lim\limits_{m\to \infty }\alpha _m=1$ where $\alpha _m=\max\limits_{k\in {\mathbb N}}|c_k^{(m})|,\, m\in \mathbb N$. }

{\bf Corollary 3.4.} {\it Let the conditions of the corollary 3.3 hold. Then there is the number  $m_0\in \mathbb N$ such that the maximum $\max\limits_{k\in {\mathbb N}}|c_{k}^{(n)}|$ is reached on the only one number  $k=i_n,\, n\in\mathbb N$ for any  $n\geq m_0$.}

{\bf Lemma 3.5}. {\it Let ${\cal F}$, ${\cal E}$ be a pair of ONB in the space $E$. 
Let $\{ m_j\}$ be a sequence of natural numbers such that $\max\limits_{i\in {\mathbb N}}|c_{i,j}|=c_{m_j,j}$ for every $j\in \mathbb N$. Then the condition 
\begin{equation}\label{up}
\sum\limits_{j=1}^{\infty}\sum\limits_{i\neq m_j}|c_{i,j}|<+\infty 
\end{equation}
is equivalent to the condition (\ref{crit}).}

{\bf Proof}. Let the condition (\ref{crit}) be hold. Let $l_j=\sum\limits_{i=1}^{\infty}|c_{i,j}|$.  Then the series  $\sum\limits_{j=1}^{\infty}\ln (l_j)$ converges. Hence the series $\sum\limits_{j=1}^{\infty} (l_j-1)$ converges and $\lim\limits_{j\to \infty }(l_j-1)=0$. 
According to lemma 3.2 the inequalities $l_j\geq \alpha _j+{\sqrt {1-\alpha _j^2}}\geq 1$ 
hold for any $j\in \mathbb N$ where $\alpha _j=\max\limits_{i}|c_{i,j}|$. 
Hence 
$l_j\geq 1-\beta _j+{\sqrt {2\beta _j-\beta _j^2}}$ for any $j\in \mathbb N$ where  $\beta _j=1-\alpha _j\in [0,1)$. 
Then $\lim\limits_{j\to \infty }\beta _j=0$ according to corollary 3.3. 
Therefore, $l_j-1 \sim {\sqrt {2\beta _j}}$ as $j\to \infty $ and the series
\begin{equation}\label{sqrt}
\sum\limits_{j=1}^{\infty }{\sqrt {\beta _j}} 
\end{equation} 
converges as well as the series  
$\sum\limits_{j=1}^{\infty }(l_j-1)$. Since  $\sum\limits_{i\neq m_j}|c_{ij}|=l_j-1+\beta _j$ for any $j\in \mathbb N$ then the condition (\ref{up}) is the consequence of the convergence of the series  $\sum\limits_{j=1}^{\infty} (l_j-1)$ and (\ref{sqrt}).

Let the condition
(\ref{up}) be hold. Then  $\lim\limits_{j\to \infty }\gamma _j=0$ where $\gamma _j=\sum\limits_{i\neq m_j}|c_{ij}|=l_j-\alpha _j,\, j\in \mathbb N$.
Since $l_j=\alpha _j+\gamma _j\leq 1+\gamma _j$ then the convergence of the series $\sum\limits_{j=1}^{\infty }\gamma _j$ is the consequence of the condition (\ref{up}). Therefore the series $\sum\limits_{j=1}^{\infty }(l _j-1)$ converges. Hence the condition (\ref{crit}) holds. \hfill$\Box $

{\bf Theorem 3.6.}
{\it The conditions (\ref{crit}) and  (\ref{crit'}) are equivalent to each other.
}

{\bf Proof}. Let the condition  (\ref{crit}) be hold. Hence $\lim\limits_{j\to \infty }l_j=1$. Therefore  $\lim\limits_{j\to \infty }\alpha _j=1$ according to the corollary 3.3. Hence there is the number $j_0\in \mathbb{N}$ such that $\alpha _j>{1\over {\sqrt 2}}$ for any $j>j_0$. 

There is the sequence of natural numbers  $m_j,\, j\in \mathbb N$ such that  $\max\limits_{k}|c_{k,j}|=c_{m_j,j}=\alpha _j$. Moreover, the number $m_j$ is uniquely defined for $j>j_0$. On the other hand there is the sequence of natural numbers 
$M_k,\, k\in \mathbb N$ such that $\max\limits_{j}|c_{k,j}|=c_{k,M_k}$. Since $|c_{m_j,j}|>{1\over {\sqrt 2}}$ for all $j>j_0$ then $\max\limits_{i}|c_{m_j,i}|=|c_{m_j,j}|$ for all $j>j_0$. I.e. the maximal element $|c_{m_j,j}|$ of $j$-th column is the maximal element of $m_j$-th row in matrix $\| C\|$ for any $j>j_0$.
Hence $M_{m_j}=j,\ \forall \ j>j_0$. 

Thus $\max\limits_{i\in {\mathbb N}}|c_{m_j,i}|=|c_{m_j,M_{m_j}}|=|c_{m_j,j}|$ for any $j>j_0$.

According to lemma $3.5$ the condition  (\ref{crit}) implies that $\sum\limits_{j\in {\mathbb N}}\sum\limits_{k\in {\mathbb N}\backslash \{ m_j\}}|c_{k,j}|<+\infty$. Therefore, \begin{equation}\label{usl}
\sum\limits_{j=1}^{j_0}\sum\limits_{k\in {\mathbb N}}|c_{k,j}|+\sum\limits_{j=j_0+1}^{+\infty }\sum\limits_{k\in {\mathbb N}\backslash \{ m_j\}}|c_{k,j}|<+\infty .
\end{equation}

Let ${\mathbb N}_1$ be the set of values of the sequence $m_j,\, j>j_0$. Let ${\mathbb N}_0$ be a set ${\mathbb N}\backslash {\mathbb N}_1$. 
In the condition (\ref{usl}) we can rearrange the order of summation of the series of non-negative terms: 
$$
\sum\limits_{j=1}^{j_0}\sum\limits_{k\in {\mathbb N}}|c_{k,j}|+\sum\limits_{j=j_0+1}^{+\infty }\sum\limits_{k\in {\mathbb N}\backslash \{ m_j\}}|c_{k,j}|=
$$
\begin{equation}\label{uslo}
=\left( \sum\limits_{k\in {\mathbb N}_0}\sum\limits_{j=1}^{j_0}+\sum\limits_{k\in {\mathbb N}_1}\sum\limits_{j=1}^{j_0}+\sum\limits_{k\in {\mathbb N}_0}\sum\limits_{j>j_0}+\sum\limits_{k\in {\mathbb N}_1}\sum\limits_{j>j_0;\, j\neq M_k}\right) |c_{k,j}|<+\infty .
\end{equation} 
In the last equality we use the following presentation of the set of summation indexes
$$\{ \{ (k,j),\ k\in {\mathbb N}\backslash \{m_j\}\},\ j>j_0 \}=\{ (k,j),\ j>j_0,\, k\in {\mathbb N}_0\}\bigcup \{ \{ (k,j),\ j>j_0,\, j\neq M_k\},\ k\in {\mathbb N}_1\}.$$
Therefore according to (\ref{uslo}) we obtain the following condition  
$$
\sum\limits_{k\in {\mathbb N}}\sum\limits_{j\in {\mathbb N}\backslash \{ M_k\}}|c_{k,j}|\leq \sum\limits_{k\in {\mathbb N}_0}\sum\limits_{j\in {\mathbb N}}|c_{k,j}|+\sum\limits_{k\in {\mathbb N}_1}\sum\limits_{j\in {\mathbb N}\backslash \{ M_k\}}|c_{k,j}|<+\infty .
$$ 
Since $\sum\limits_{k\in {\mathbb N}}\sum\limits_{j\in {\mathbb N}\backslash \{ M_k\}}|c_{k,j}|<+\infty $ then according to the lemma 3.5 the condition (\ref{crit'}) holds. 
If we swap the bases $\cal E$ and $\cal F$ then we obtain that the condition (\ref{crit'}) implies the condition (\ref{crit}). \hfill$\Box$

{\bf Corollary 3.7}. {\it Let $\|C\|$ be the transition matrix from one ONB to another. Then the product $\prod\limits_{j=1}^{\infty }\| c_{\cdot ,j}\|_{\ell _1}$ converges if and only if the product $\prod\limits_{i=1}^{\infty }\| c_{i,\cdot }\|_{\ell _1}$ converges.}

ONB $\cal E$ is called near to the ONB $\cal F$ if  $\cal E$ and $\cal F$ satisfy the condition (\ref{crit}).
If the condition (\ref{crit}) is not satisfy for two ONB $\cal E$ and $\cal F$ then  ONB $\cal E$ is called distant from the ONB $\cal F$.

In the section 4 we show that if ONB $\cal E$ is near to ONB $\cal F$ then the measures  $\lambda _{\cal E}$ and $\lambda _{\cal F}$ coincide. In the section 5 we show that if ONB $\cal E$ is distant from ONB $\cal F$ then measures  $\lambda _{\cal E}$ and $\lambda _{\cal F}$ are defined on the different rings such that both measures $\lambda _{\cal E}$ and $\lambda _{\cal F}$ take zero values on an arbitrary set from the ring ${\cal R}_{\cal E}\bigcap {\cal R}_{\cal F}$. The results of sections 4, 5 are obtained by the technique analysis of mutual position of rectangles with the edges collinear to the vectors of different ONB. 
The results of sections 4, 5 gives the proof of the existence and describes the properties of shift and rotation-invariant measures on the Hilbert space. This  proof do not use the choice axiom as opposed to the approach of the article \cite{s16}.

\section{Ring ${\cal R}_{\cal E}\bigcap {\cal R}_{\cal F}$ in the case of nearness of bases $\cal E$ and $\cal F$} 

Let us prove that if the condition (\ref{crit}) holds  (as well as the equivalent condition (\ref{crit'})) then the rectangles $\Pi $ and $Q={\bf U}(\Pi )$ belong to the ring ${\cal R}_{\cal E}$ and $\lambda _{\cal E}(Q)=\lambda _{\cal E}(\Pi ).$ At the first step to this goal we apply the orthogonal mappings  
$\bf V$ of the space $E$ which only change the order of vectors in the basis $\cal E$. Also we use the following property of the measure $\lambda _{\cal E}$ to be invariant with respect to permutation of vectors of basis $\cal E$.

{\bf Theorem 4.1}. {\it Let the orthogonal transformation $\bf V$ of the space $E$ changes the order of vectors of the basis
$\cal E$ only and ${\cal E}'={\bf V}({\cal E})$. Let the orthogonal transformation  $\bf W$ of the space  $E$ changes the direction of some vectors of the basis 
 $\cal E$ only and ${\cal E}''={\bf W}({\cal E})$. Then 
${\cal R}_{{\cal E}''}={\cal R}_{{\cal E}'}={\cal R}_{\cal E}$ and $\lambda _{{\cal E}''}(A)=\lambda _{{\cal E}'}(A)=\lambda _{\cal E}(A)$ for all $A\in {\cal R}_{\cal E}$.}

{\bf Proof}. The collections of absolutely measurable rectangles  ${\cal K}_{\cal E}$ and ${\cal K}_{{\cal E}'}$ coincide with each other. In fact, for any rectangle 
$\Pi \in {\cal K}_{\cal E}$ its edges are collinear to the vectors of ONB ${\cal E}'$ and the product of the lengths of edges converges unconditional. Therefore $\Pi \in {\cal K}_{{\cal E}'}$. The opposite is also true. Moreover, 
  $\lambda _{{\cal E}'}(\Pi )=\lambda _{\cal E}(\Pi )$ for any $\Pi \in {\cal K}_{\cal E}$ according to unconditional measurability. 
Since  ${\cal K}_{{\cal E}'}={\cal K}_{\cal E}$ then $r_{\cal E}=r_{{\cal E}'}$. Since the finitely additive function $\lambda _{{\cal E}'}$ coincides with  $\lambda _{\cal E}$ on the collection of sets
 ${\cal K}_{{\cal E}'}={\cal K}_{\cal E}$ then the additive functions
$\lambda _{{\cal E}'}$ and $\lambda _{\cal E}$ has the only one additive extension onto the ring 
$r_{\cal E}=r_{{\cal E}'}$.
Therefor the completion of the measures $\lambda _{{\cal E}'}$ and  $\lambda _{\cal E}$ coincides with each other. Thus we prove the statement on the transformation   $\bf V$. The statement on the transformation $\bf W$ can be proved analogously. \hfill$\Box$

\medskip

Let ${\cal E}=\{ e_j\}$ and  ${\cal F}=\{ f_k\}$ be a pair of ONB in the space $E$. Let us define subspaces 
$E_n={\rm {span}}(e_1,\ldots, e_n)$ and $E^n=E_n^{\bot }$; $F_n={\rm {span}}(f_1, \ldots ,f_n)$ and $F^n=F_n^{\bot }$ for any $n\in {\mathbb N}$.

\medskip

{\bf Theorem 4.2.}
{\it Let the condition (\ref{crit}) holds for the pair of ONB   $\cal E$ and $\cal F$. Then there is the ONB ${\cal F}'$ such that ${\cal F}'$ can be obtained  by means of changing of numbering of vectors  ONB $\cal F$  and satisfying the condition
\begin{equation}\label{dia}
\exists \, m_0\in {\mathbb N}:\ \forall \ j>m_0\;\; \ |c'_{j,j}|=\max\limits_{k\in {\mathbb N}}|c_{j,k}'|=\max\limits_{i\in {\mathbb N}}|c_{i,j}'|, 
\end{equation}
where $c'_{k,j}=(f_k',e_j),\, k,j\in {\mathbb N}.$
}

\medskip

{\bf Proof}. According to condition (\ref{crit}) we have $\lim\limits_{j\to \infty }\alpha _j=1$. Therefor the set $M'=\{ j\in {\mathbb N}:\ \alpha _j\leq {1\over {{\sqrt 2}}} \}$ is finite. Since the condition (\ref{crit}) implies the condition (\ref{crit'}) then the set  $M''=\{ k\in {\mathbb N} :\ \max\limits_{j\in {\mathbb N}}|c_{k,j}|\leq {1\over {{\sqrt 2}}} \}$ is finite analogously. Let  $m',m''$ be a numbers of elements in finite sets  $M',M''$ respectively.  Let us prove that  $m'=m''$.

Let us assume that $m'>m''$ (the case $m'<m''$ can be considered analogously).
Then  $|c_{k,j}|\leq {1\over {\sqrt 2}}$ for any $ k\in M'',\, j\in M'$.

We will done finite number of permutation of rows and columns of the matrix $\| c_{ij}\|$. For each permutation of two rows (of two columns) we done the permutation of corresponding vectors in the basis $\cal F$ (in the basis $\cal E$).

Step 1. Let us permute the vectors of ONB $\cal F$ with the numbers from the set  $M''$ onto first $m''$ positions. The natural order of numbers in the set   $M''$ and  in its complement are preserved. Analogously, let us permute the vectors of ONB $\cal E$ with the numbers from the set  $M'$ onto first $m'$ positions. The natural order of numbers in the set   $M'$ and  in its complement are preserved.

After this permutation of bases $\cal E$ and $\cal F$ we obtain the matrix of transition with the following properties. Each row with number greater than $m''$ contains the only one element with modulus greater than  ${1\over {\sqrt 2}}$. According to the permutation in step 1 this element belongs to the column with the number greater than $m'$. Conversely,  each column with number greater than $m'$ contains the only one element with modulus greater than  ${1\over {\sqrt 2}}$.  This element belongs to the row with the number greater than  $m''$. 

Hence there is the permutation of columns with the numbers  $m'+1,m'+2,\ldots \,$ such that each row with a number  $k>m''$ contains the only one element with modulus greater than ${1\over {\sqrt 2}}$ and this element belong to the column with the number $k+m'-m''$.

Step 2. We transform the ONB $\cal F$ by the following rule. We change the vector $f_k,\, k>m''$, onto the vector  $-f_k$ under the condition $c_{k,k+m'-m''}<0$. In opposite case we remain $f_k$. 

The ring of subsets
${\cal R}_{\cal F}$ (and ${\cal R}_{\cal E}$) and the measure $\lambda _{\cal F}$ (and $\lambda _{\cal E}$) so not change under the transformation in the steps 1) and 2) according to the theorem 4.1. The matrix $\| c_{i,j}\|$ under the above transformation satisfy the conditions  $c_{k,k+m'-m''}>{1\over {\sqrt 2}}\;\; \forall \ k>m''$.

In the proof of the lemma 3.5 we should prove that the condition (\ref{crit}) implies  estimates  
$$\sum\limits_{j>m''}(1-c_{j,j+m'-m''})<\sum\limits_{k=1}^{\infty }(1-|c_{k,M_k}|)<\sum\limits_{k=1}^{\infty }{\sqrt{1-|c_{k,M_k}|}}<+\infty.$$ 
Therefore, there is the number $N$ such that 
$\sum\limits_{j=N+1}^{\infty}{\sqrt {1-c_{j,j+m'-m''}}}<{1\over 2}$ and\\ $\sum\limits_{k=N+1}^{\infty}(l_k^T-|c_{k,k+m'-m''}|)<{1\over 2}\, $  (here $l_k^T=\sum\limits_{j=1}^{\infty }|c_{kj}|$).

Let $E^{N+m'-m''}$ be the subspace of the space $E$ such that orthonormal system of vectors $\{ e_{N+m'-m''+1},\ldots \}$ forms the ONB in the space $E^{N+m'-m''}$. Let ${\bf P}_{E_a}$ be an orthogonal projector in the space $E$ onto a subspace $E_a$ of the space $E$. Let us consider the system of vectors 
$\tilde f_{N+1},\tilde f_{N+2},\ldots$. Here  $\tilde f_k={\bf P}_{E^{N+m'-m''}}f_k$ for any $k>N$. Therefore $\|f_k-\tilde f_k\|_E\leq l^T_k-|c_{k,k+m'-m''}|$ and
$$
\|e_{k+m'-m''}-\tilde f_k\|\leq \Bigl\| \sum\limits_{i\neq k+m'-m''}c_{i}e_i+(1-c_{k+m'-m''})e_{k+m'-m''}\Bigr\| 
={\sqrt {2(1-c_{k,k+m'-m''})}}<1.
$$

Hence the system of vectors $\{ \tilde f_k,\, k>N\}$ of the space $E^{N+m'-m''}$ is the perturbation of ONB $\{ e_k,\, k>N+m'-m''\}$ which is small in the following sense $\sum\limits_{j=N+1}^{\infty }\| e_{j+m'-m''}-\tilde f_j\|_E<1$. Therefore the system of vectors $\{ \tilde f_k,\, k>N\}$ is the Riesz  basis in the space $E^{N+m'-m''}$ (see \cite{Vilenkin}, chapter 1.6).  
The system of vectors $\{ f_1,\ldots,f_N,\tilde f_{N+1},\ldots\}$ is the Riesz basis in the space  $E$ since it is nearby to the ONB  ${\cal F}=\{f _j\}$ in the following sense
$\sum\limits_{j=N+1}^{\infty }\| f_{j}-\tilde f_j\|_E\leq \sum\limits_{j=N+1}^{\infty}(l^T_j-|c_{j,j+m'-m''}|)<1$. Since system of vectors $\{\tilde f_{N+1},\ldots \}$ forms Riesz basis in the subspace  $E^{N+m'-m''}$ of codimension  $N+m'-m''$ then the system of vectors $\{ \hat f_1,\ldots,\hat f_N,\tilde f_{N+1},\ldots\}$  forms Riesz basis in the space $E$ (here $\hat f_i={\bf P}_{E_{N+m'-m''}}f_i,\, i=1,\ldots,N$). Therefore the system of vectors $\{ \hat f_1,\ldots,\hat f_N\}$ form the basis in the space $E_{N+m'-m''}$. It is impossible in the case $m'>m''$. The contradiction proves that $m''=m'$. \hfill$\Box$

{\bf Remark 4.3.} The conditions (\ref{crit}), (\ref{crit'}) and (\ref{up}) are invariant with respect to changing of numbering of vectors of bases ${\cal E}, \ {\cal F}$.

{\bf Remark 4.4.} 
The theory of determinants of linear operators \cite{FB} forms different approaches to definition of a determinant  and to study a conditions of its existence.
The Poincare theorem gives the condition on the infinite matrix of a linear operator in some basis sufficient to the existence of determinant. Poincare theorem  (see \cite{ArSad}, p. 400) states that the following two conditions A) and B) are sufficient for the existence of the determinant of infinite matrix  $\|C\|$ (the determinant of the infinite matrix  $\|C\|$ is defined as the limit of  $n$-th order main angular minor of matrix $\| C\|$ as $n\to \infty$). 
Here A) is the condition of unconditional convergence of the products of diagonal matrix elements; B) is the condition of absolute convergence of the series of non-diagonal elements of matrix $\| C\|$.
If the matrix  $\| C\|$ is orthogonal then the condition A) is the consequence of the condition B).
In this case the condition (\ref{up}) on the pair of ONB $\cal E$ and $\cal F$ is equivalent to the condition B) of the theorem 25 \cite{ArSad}  for the matrix which is connected with matrix  $\| C\|=\| (f_k,e_j)\|$ by means of permutations of rows and columns. Thus the deformation of the measure under the action of linear mapping connected with the determinant of the mapping (see also \cite{NOMA}).

\bigskip

Let us introduce some notations.
Let ${\cal E},\, {\cal F}$ be a pair of ONB.
Finite-dimensional subspaces $E_n={\rm span }(e_1,\ldots , e_n)$, $F_n={\rm span }(f_1,\ldots , f_n)$ and its orthogonal completion  $E^n=(E_n)^{\bot}$ и $F^n=(F_n)^{\bot}$ are defined for any $n\in \mathbb N$. 
The operators ${\bf P}_{E_n}$ and ${\bf P}_{E^n}$ in the space $E$ are the operators of orthogonal projections onto the subspaces $E_n$ and $E^n$ respectively. 
For any $n\in {\mathbb N}$ the matrix $ C_n
=\| (e_i,f_j)\|, \ i,j\in \{1, \ldots , n \}$, is the matrix of orthogonal projectors ${\bf P}_{F_n,E_n}:\; F_{n}\to E_n$ from the subspace $F_n$ into the subspace $E_n$ in pair of bases $\{ f_1,...,f_n\}$ and $\{ e_1,...,e_n\}$. Let $\lambda _n$ be the Lebesgue measure in $n$-dimensional euclidean space.

{\bf Lemma 4.5.} {\it Let  ${\cal E}, \ {\cal F}$ be a pair of bases which satisfy the conditions  (\ref{crit}) and (\ref{dia}). Then for any $\epsilon >0$ there is a number $N_{\epsilon}\in {\mathbb N}$ such that  ${\rm Tr}(C _n^TC_n)\geq n-\epsilon$ for any $n\geq N_{\epsilon}
$.}

To prove lemma 4.5 we firstly obtain some asymptotic estimates for the spectrum of the matrix 
$C_n
=\| (e_i,f_j)\|, \ i,j\in \{1, \ldots , n \}$.

The norms of projections of vectors $f_j,\, j=1,\ldots,n$ onto a subspace $E_n$ has the following expression 
$\|{\bf P}_{E_n}f_j\| ^2=\sum\limits_{k=1}^nc_{kj}^2$.
Therefore
\begin{equation}\label{>m}
1\geq \|{\bf P}_{E_n}f_j\| ^2\geq c_{jj}^2=\alpha _j^2 \;\; \forall \ j\in \{1, \ldots , n\}.
\end{equation}

Let the number $\epsilon >0$ be fixed. Since the series $\sum\limits_{k=1}^{\infty }(1-\alpha _k)$ converges (lemma 3.5) then there is  $m\in \mathbb N$ such that 
\begin{equation}\label{us1}
\sum\limits_{k=m+1}^{n}(1-\alpha _k)\leq {{\epsilon }\over 2}
\end{equation}
for any $n\geq m$. Then according to (\ref{>m}) and  (\ref{us1}) we have the estimate  
$$
\sum\limits_{k=m+1}^{n }(1-\| {\bf P}_{E_n}f_k)\| ) \leq \sum\limits_{k=m+1}^{\infty }(1-\alpha _k)\leq {{\epsilon }\over 2}\quad \forall \ n>m.
$$

Since the sequence of operators $\{ {\bf P}_{E_n}\}$ converges to unit operator in the strong operator topology then $\lim\limits_{n\to \infty }\| {\bf P}_{E_n}f_j\|=1$ for any $j\in \{ 1,\ldots ,m\}$. Therefore there is a number $N_{\epsilon }>m$ such that the inequality  
\begin{equation}\label{us2}
\sum\limits_{j=1}^m(1- \|{\bf P}_{E_n}f_j\| )<{{\epsilon }\over 2}
\end{equation}
holds for any  $n\geq N_{\epsilon }$.
Thus for any $\epsilon >0$ there is a numbers $m_{\epsilon }\in \mathbb N$ and $N_{\epsilon}>m_{\epsilon}$ such that the condition 
\begin{equation}\label{us3}
\sum\limits_{j=1}^{n}(1- \|{\bf P}_{E_n}f_j\| )=\sum\limits_{j=1}^{m_{\epsilon }}(1- \|{\bf P}_{E_n}f_j\| )+\sum\limits_{j=m_{\epsilon }+1}^{n}(1- \|{\bf P}_{E_n}f_j\| )<{{\epsilon }} 
\end{equation}
hods for all  $n\geq N_{\epsilon}$.

The columns of the matrix $\| C_n\|$ are the coordinate columns of vectors ${\bf P}_{E_n}f_1,\ldots,{\bf P}_{E_n}f_n$ with respect to ONB $\{ e_1, \ldots, e_n\}$ of the space $E_n$. Hence the equality  $(C_n^TC_n)_{jj}=\|{\bf P}_{E_n}f_j\|^2$ holds for any $j\in \{1, \ldots , n\}$. Therefore,
$$
{\rm Tr}(C^T_nC_n)=\sum\limits_{j=1}^n\| {\bf P}_{E_n}f\|_j^2=\sum\limits_{j=1}^n(1-(1-\| {\bf P}_{E_n}f\|_j))^2\geq n-2\sum\limits_{j=1}^n(1-\| {\bf P}_{E_n}f\|_j)
$$ 
for any $n \in \mathbb N $.  Thus according to  (\ref{us3}) there is a number  $N_{\epsilon}\in {\mathbb N}$ such that   
${\rm Tr}(C^T_nC_n) \geq n-2\epsilon $ for any  $n\geq N_{\epsilon}$. \hfill$\Box$

{\bf Corollary 4.6.} {\it Let the assumption of lemma 4.5 be hold. Then for any $\epsilon \in (0,{1\over 2})$ there is a number $N_{\epsilon}\in {\mathbb N}$ such that ${\rm det}(C_n^TC_n)\geq 1-2\epsilon$ for any $n\geq N_{\epsilon}$.}

The matrix $C_n^TC_n$ is positive   $n\times n$ matrix. Hence it has the ONB of eigenvectors and the collection of $n$ positive eigenvalues $\mu _1,\ldots,\mu _n$. Since $C_n$ is the matrix of the operator of orthogonal projection  ${\bf P}_{F_n,E_n}:\ F_n\to E_n$ (and $C_n^T$ is the matrix of operator of orthogonal projection ${\bf P}_{E_n,F_n}$ from the subspace $E_n$ onto the subspace $F_n$) then the eigenvalues of the matrix  $C^T_nC_n$ are no greater than 1 since $\|{\bf P}_{F_n,E_n}x\|_E\leq \|x\|_F \ \forall \ x\in F_n$ and $\|{\bf P}_{E_n,F_n}y\|_F\leq \|y\|_E \ \forall \ y\in E_n$. Thus $0<\mu _k\leq 1\ \forall \ k\in \{1,\ldots , n\}$.  

Let us fix some $\epsilon \in (0,{1\over 2})$. Then according to lemma 4.5 there is the number  $N_{\epsilon}\in {\mathbb N}$ such that the inequality 
${\rm Tr}(C^T_nC_n)=\sum\limits_{k=1}^n(1-(1-\mu _k))=n-\sum\limits_{k=1}^n(1-\mu _k)\geq n-\epsilon$ holds for any $n\geq N_{\epsilon}$. Hence $\sum\limits_{k=1}^nt_k<\epsilon$ where $t_k=1-\mu _k\in [0,\epsilon ]$. Therefore ${\rm det}(C_n^TC_n)=\prod \limits_{k=1}^n\mu _k=\prod\limits_{k=1}^n(1-t_k)
\geq 1+\sum\limits_{k=1}^n\ln (1-t_k)$. According to Lagrange theorem we have $\ln (1-t)>-{t\over {1-\epsilon }}\;\; \forall \ t\in [0,\epsilon]$. Therefore the inequality ${\rm det}(C_n^TC_n)\geq 1-2\epsilon$ holds for any  $n\geq N_{\epsilon}$. 
\hfill $\Box$

{\bf Lemma 4.7.} {\it  Let  ${\cal E}, \ {\cal F}$ be a pair of bases which satisfy the conditions  (\ref{crit}) and (\ref{dia}). Then $\lim\limits_{n\to \infty}\lambda _n({\bf P}_{F_n,E_n}(Q_{n}))=1.$}

In fact, $\lambda _n({\bf P}_{F_n,E_n}(Q_{n}))=|{\rm det}(C_n)|$. According to corollary 4.6 for any $\epsilon \in (0,{1\over 2})$ there is the number $N_{\epsilon}\in {\mathbb N}$ such that the condition ${\rm det}(C_n^TC_n)=({\rm det}(C_n))^2\geq 1-\epsilon$ holds for any $n\geq N_{\epsilon}
$. Hence the estimates  $||{\rm det}(C_n)|-1|=1-|{\rm det}(C_n)|$ hold for any $n\geq N_{\epsilon}
$. Thus the statement of lemma 4.7 is proved. \hfill$\Box$

{\bf Remark 4.8.} The proof of the lemma 4.7 is based on the existence of the determinant ${\rm det}(C)=\lim\limits_{n\to \infty }{\rm det}(C_n)$ of the matrix $C=\| (f_k,e_j)\|$ (see \cite{ArSad, FB, GK65, Kuk}).

{\bf Theorem 4.9}. {\it Let the condition  (\ref{crit}) be hold. If  $Q\in {\cal K}_{\cal F}$ then $Q\in {\cal R}_{\cal E}$ and $\lambda _{\cal E}(Q)=\lambda _{\cal F}(Q)$. On the contrary, if $P\in {\cal K}_{\cal E}$ then $P\in {\cal R}_{\cal F}$ and $\lambda _{\cal F}(P)=\lambda _{\cal E}(P)$.}

{\bf Proof}. According to the theorem 4.2 we can count that the condition (\ref{dia}) holds. (In opposite case we can change the numbering of basis vectors to obtain (\ref{dia}) but the rings ${\cal R}_{\cal E}$, ${\cal R}_{\cal F}$ and measures $\lambda _{\cal E},\, \lambda _{\cal F}$ do not change). 
If we prove the first statement then the second one follows from the first statement and the theorem 3.6.
To prove the first statement of the theorem 4.9 
it is sufficient to prove the statement for a rectangle $Q\in {\cal K}_{\cal F}$ such that $\lambda _{\cal F}(Q)>0$.

In fact, let $Q\in {\cal K}_{\cal F}$ and $\lambda _{\cal F}(Q)>0$. Then for any $\epsilon >0$ the rectangle $Q\in {\cal K}_{\cal F}$ can be inscribed into the measurable rectangle  $Q'\in {\cal K}_{\cal F}$ with positive measure $\lambda _{\cal F}(Q')<\epsilon $. Hence $Q'\in {\cal R}_{\cal E}$ and the external measure $\bar\lambda _{\cal E }$ admits estimates $\bar\lambda _{\cal E }(Q)\leq \lambda _{\cal E}(Q')=\lambda _{\cal F}(Q')<\epsilon $. Since $\epsilon >0$ is arbitrary then $\lambda _{\cal E }(Q)=0=\lambda _{\cal F }(Q)$.

Let us show that if $Q\in {\cal K}_{\cal F}$ and $\lambda _{\cal F}(Q)>0$ then for any $\epsilon >0$ there are

1) a set $S \in {\cal R}_{\cal E}$ such that  $Q\subset S$ and $\lambda _{\cal E}(S)<(1+\epsilon )\lambda _{\cal F}(Q)$ (the estimate from above); 

2) a set  $s \in {\cal R}_{\cal E}$ such that $s\subset Q$ and $\lambda _{\cal E}(s)>(1-\epsilon )\lambda _{\cal F}(Q)$ (the estimate from below).

{\bf I}. The estimate from above. Let $Q\in {\cal K}_{\cal F}$ and $\lambda _{\cal F }(Q)>0$.
Let $\{ d_k\}$ be a sequence of lengths of the edges of the rectangle $Q$. Then $\sum\limits_{k=1}^{\infty }\max \{0, d_k-1\}<\infty $ since rectangle $Q$ is measurable. Moreover, the series $\sum\limits_{k=1}^{\infty }(d_k-1)$ converges absolutely since   $\lambda _{\cal F}(Q)>0$. Hence $0<d_0\leq 1\leq D_0<+\infty $ where $d_0=\inf \{ d_k\},\, D_0=\sup\{ d_k\}$. Since the measures $\lambda _{\cal E},\, \lambda _{\cal F}$ are invariant with respect to shift then we can count that the rectangle $Q$ is centered and it can be parametrized by the following way
$$
Q=\{ x\in E:\ (x,f_k)\in [-{1\over 2}d_k,{1\over 2}d_k],\, k\in {\mathbb N}\} .
$$
Since the length of the rectangle $Q$  projection onto the axis $Oe_j,\, j=1,2,\ldots $ is equal to $b_j=\sum\limits_{k=1}^{\infty }d_k|c_{k,j}|,\ j\in {\mathbb N}$ then the rectangle $Q$ can be inscribed into the rectangle $\Pi $ such that the edges of the rectangle $\Pi$ are collinear to the vectors of ONB $\cal E$ and the lengths of edges form the sequence $\{ b_j\}$.

Let us prove that the rectangle $\Pi \in {\cal K}({\cal E})$ and $\lambda _{\cal E}(\Pi )>0$. It is sufficient to prove that the series $\sum\limits_{j=1}^{\infty }|b_j-1|$ converges. For every $j\in \mathbb N$ the equalities $b_j=\sum\limits_{k=1}^{\infty }d_k|c_{k,j}|=\sum\limits_{k=1}^{\infty }|c_{k,j}|+\sum\limits_{k=1}^{\infty }(d_k-1)|c_{k,j}|=l_j+\sum\limits_{k=1}^{\infty }(d_k-1)|c_{k,j}|.$ hold. Hence we have the estimates
\begin{equation}\label{Pi-est}
\sum\limits_{j=1}^{\infty }|b_j-1|\leq \sum\limits_{j=1}^{\infty }(l_j-1)+\sum\limits_{j=1}^{\infty }\sum\limits_{k=1}^{\infty }|d_k-1)||c_{k,j}|=\sum\limits_{j=1}^{\infty }(l_j-1)+\sum\limits_{k=1}^{\infty }|d_k-1)|\hat l_k<+\infty .
\end{equation}
In fact, $l_j=\sum\limits_{k=1}^{\infty }|c_{k,j}|>1 \ \forall \ j$ and the series $\sum\limits_{j=1}^{\infty }(l_j-1)$ converges according to condition (\ref{crit}). Analogously, $\hat l_k=\sum\limits_{j=1}^{\infty }|c_{k,j}|>1 \ \forall \ k$ and the series $\sum\limits_{k=1}^{\infty }(\hat l_k-1)$ converges according to the theorem 3.6. Hence the sequence $\{\hat l_k\}$ is bounded. Therefore the estimate (\ref{Pi-est}) holds. 

Let us construct the centered rectangle ${\tilde \Pi }\in {\cal K}_{\cal E}$ such that the lengths of the edges of this rectangle form the sequence $\{ {\tilde b}_j\},\ \tilde b_j=\max \{ 1,b_j\} ,\, j\in \mathbb N$. Then the rectangle $\tilde \Pi $ is absolutely measurable according to the estimate (\ref{Pi-est}).

Since the rectangle  $Q\in{\cal K}_{\cal F}$ is measurable then  the equality $\lambda _{\cal F}(Q)=\lim\limits_{n\to \infty }\lambda _{n}(Q_n)$ holds. Here  
$Q_n$ is the $n$-dimensional section of the rectangle $Q$ by the hyperplane $F_n$  for each  $n\in \mathbb N$.

Let ${p}_n={\bf P}_{F_n,E_n}(Q_n)$ be the orthogonal projection of the $n$-dimensional rectangle $Q_n=F_n\bigcap Q$  from the subspace $F_n$ onto the subspace $E_n$. The matrix $\| c^{(n})_{i,j}\|,\, i,j\in {\overline {1,n}}$ is the Jacobi matrix of the linear mapping of orthogonal projection  ${\bf P}_{F_n,E_n}$ in the bases $f_1,\ldots,f_n$ and $e_1,\ldots,e_n$ in the subspaces  $F_n$ and $E_n$. Therefore $\lambda _n(p_n)=|\det (\| c^{(n)}_{i,j}\|)|\lambda _n(Q_n)$.

Let $P_n={\bf P}_{E_n}(Q)$ be the orthogonal projection onto the subspace $E_n$ of the unit rectangle $Q$. Let $\tilde \Pi _n$ be the projection of the rectangle $\tilde \Pi$ onto  $n$-dimensional hyperplane $E_n$. Then  $p_n\subset P_n\subset \tilde \Pi _n$ for any $n\in \mathbb N$ and $P_n$ be the convex subset of the space $E_n$. 
Let
\begin{equation}\label{B}
\prod _{j=1}^{\infty }\tilde b_j=B<\infty .
\end{equation} 
Hence $\lambda _n(\tilde \Pi _n)=\tilde b_1 \cdots \tilde b_n\leq B\ \forall \ n \in \mathbb N$.

The rectangle $Q$ can be parametrized by the shifts of $n$-dimensional rectangle  $Q_n$ along the vectors $f_{n+1},\ldots$ by the equality  $Q=\Bigl\{ Q_n+\sum\limits_{k=n+1}^{\infty}d_kt_kf_k,\ t_k\in [-{1\over 2},{1\over 2}]  ,\, \sum\limits_{k=n+1}^{\infty }d_k^2t_k^2<\infty \Bigr\}$. Therefore the set $P_n$ admits the parametrization
\begin{equation}\label{s}
\Bigl\{ p_n+\sum\limits_{k=n+1}^{\infty }d_kt_k{\bf P}_{E_n}(f_k),\ t_k\in \Bigl[-{1\over 2},{1\over 2}\Bigr] ,\, \sum\limits_{k=n+1}^{\infty }d_k^2t_k^2<\infty \Bigr\}
=\Bigl\{ p_n+\sum\limits_{j=1}^nt_j\beta _{n,j}e_j,\ t_j\in [-{1\over 2},{1\over 2}]\Bigr\}
\end{equation} 
where $\beta _{n,j}=\sum\limits_{k=n+1}^{\infty }d_k|c_{k,j}|$.  For every $n\in \mathbb N$ we define the numbers
\begin{equation}\label{delt}
\Delta _{j,n}=\sum\limits_{k=n+1}^{\infty }|c_{k,j}|,\ j=1,\ldots ,n.
\end{equation}
Hence the set $P_n$ lies in the convex hull $S_n$ of shifts of the set $p_n$ along the axes $Oe_j$ onto the vectors $tD_0\Delta _{j,n}e_j,\, t\in [-{1\over 2},{1\over 2}],\, j=1,\ldots,n$:
$$
P_n\subset  \bigcup\limits_{(t_1,\ldots,t_n)\in [-1,1]^n}(p_n+\sum\limits_{j=1}^{n}{1\over 2}t_jD_0\Delta _{j,n}e_j)  
\equiv S_n\subset  \tilde \Pi _n.
$$

For every $n\in \mathbb N$ the set $S_n$ is the result of sequentially for $j=1,\ldots ,n$ elongations of the set $p_n$ along the axis $Oe_j$ onto the value ${1\over 2}D_0\Delta_{j,n}$ in the directions of vectors $e_j$ and $-e_j$. We prove that the set $S_n$ is Jordan measurable and obtain the estimate 
\begin{equation}\label{indu}
\lambda _n(S_{n})\leq \lambda (p_n)+\sum\limits_{j=1}^nBD_0\Delta _{n,j}
\end{equation}
for its Jordan measure by using the induction with respect to sequentially elongation.

Let $S_{n,0}=p _n,\ S_{n,k}=\bigcup\limits_{(t_1,\ldots,t_k)\in [-1,1]^k}(p_n+{1\over 2}\sum\limits_{j=1}^{k}t_jD_0\Delta _{j,n}e_j) $ for every $k=1,...,n. $ Then $S_n=S_{n,n}$. We prove that the sets $S_{n,k},\, k=1,...,n$ are Jordan measurable and obtain the estimates (\ref{indu}) by the induction with respect to indexes $k\in \{ 0,1,...,n\}$.
The set $S_{n,0}$ is Jordan measurable and $\lambda _n(S_{n,0})=\lambda (p_n).$
Let $k\in \{ 0,...,n-1\}$. Let the set $S_{n,k}$ is Jordan measurable and the estimate 
\begin{equation}\label{induc}
\lambda _n(S_{n,k})\leq \lambda )p_n+\sum\limits_{j=1}^kBD_0\Delta _{n,j}
\end{equation}
holds. Then $S_{n,k+1}=S_{n,k}\bigcup \Delta S_{n,k}$ where $$\Delta S_{n,k}=S_{n,k+1}\backslash S_{n,k}=\bigcup\limits_{t_{k+1}\in [-1,1]}(S_{n,k}+{1\over 2}t_{k+1}D_0\Delta _{n,k+1}e_{k+1}).$$
Then the set $\Delta S_{n,k}$ is Jordan measurable.
Since the orthogonal projection ${\bf P}_{\Gamma _{k+1}}$ of the set $S_{n,k}$ onto the $(n-1)$-dimensional hyperplane $\Gamma _{k+1}={\rm span}(e_1,...,e_{k},e_{k+2},...,e_n)$ belongs to the projection of the rectangle $\tilde \Pi _n$ onto this hyperplane then $\lambda _{n-1}({\bf P}_{\Gamma _{k+1}}(S_{n,k}))\leq \lambda _{n-1}({\bf P}_{\Gamma _{k+1}})(S_{n,k})\leq \tilde b_1...\tilde b_k\tilde b_{k+2}...\tilde b_n\leq B$. Since the elongation of the set $S_{n,k}$ onto the set $S_{n,k+1}$ along the axis $Oe_{k+1}$ is equal to $D_0\Delta _{n,k+1}$ then $\lambda _n(\Delta S_{n,k})\leq BD_0\Delta _{n,k+1}$. Therefore the sets $S_{n,k}$ are Jordan measurable and the estimates (\ref{induc}) are hold for any $k=0,1,...,n$.

Thus for every $n\in \mathbb N$ the set $S_n$ is measurable and the estimate (\ref{indu}) takes place.

{\bf Lemma 4.10.} {\it Let the condition  (\ref{crit}) be held. Then
$\lim\limits_{n\to \infty }\Bigl(\sum\limits_{j=1}^n\Delta _{j,n}\Bigr)=0$.
}

Let us fix a number $\epsilon >0$. Since the series  $\sum\limits_{j=1}^{\infty }(1-\alpha _j)$ converges then there is the number $m_{\epsilon }\in \mathbb N$ such that $\sum\limits_{j=1+m}^{\infty }(1-\alpha _j)<{{\epsilon}\over 2}$ for every $m\geq m_{\epsilon}$.
According to the condition (\ref{crit})  the inequality $\sum\limits_{k=1}^{\infty }|c_{k,j}|<+\infty $ holds for any $j\in \mathbb N$. Hence 
there is the number $N
>m_{\epsilon }$ such that $\sum\limits_{j=1}^{m_{\epsilon }}\sum\limits_{k=N+1}^{\infty }|c_{k,j}|<{{\epsilon}\over {2BD_0}}$.
Thus for every $n>N$ we have the estimates 
$$
\sum\limits_{j=1}^n\Delta _{j,n}=\sum\limits_{j=1}^{m_{\epsilon }}\Delta _{j,n}+\sum\limits_{j=m_{\epsilon}+1}^n\Delta _{j,n}=\sum\limits_{j=1}^{m_{\epsilon }}\sum\limits_{k=n+1}^{\infty }|c _{k,j}|+\sum\limits_{j=m_{\epsilon}+1}^{n}\sum\limits_{k=n+1}^{\infty }|c _{k,j}<
$$ 
$$<\sum\limits_{j=1}^{m_{\epsilon }}\sum\limits_{k=N+1}^{\infty }|c_{k,j}|+\sum\limits_{j=m_{\epsilon}+1}^{n}(1-\alpha _j)<\epsilon .\eqno{\square}
$$

{\bf Corollary 4.11.} {\it Let the condition  (\ref{crit}) be held. Then for every $n\in \mathbb N$ there is the measurable set $S_n\subset E_n$, such that $P_n\subset S_n$ and the sequence $\{ S_n\}$ satisfy the condition $\lim\limits_{n\to \infty}(S_n)=\lambda _{\cal F}(Q)$.}

The statement of corollary 4.11 is the consequence of the lemmas 4.7, 4.10 and the estimates (\ref{indu}). \hfill$\Box$

\medskip 

Let $B_j$ be the projection of the rectangle $Q$ onto the axis $Oe_j$. Then $B_j$ is the segment with the length $b_j=\sum\limits_{k=1}^{\infty }d_k|(e_j,f_k)|$). 
Since $P_n$ is the projection of the rectangle $Q$ onto the hyperplane  $E_n$ and $P_n\subset S_n$, 
then
$$
Q\subset S_n\times {B}_{n+1}\times {B}_{n+2}\times \cdots\, \quad \forall \ n\in {\mathbb N}.
$$
According to (\ref{Pi-est})
the series $\sum\limits _{k=1}^{\infty }(b_k-1)$ converges absolutely.

Let us fix a number $\epsilon >0$. Then there is the number $n_0\in \mathbb N$ such that, at first, $\prod\limits_{j=n+1}^{\infty }b_j\in [1,1+{{\epsilon }\over 2})$ for any $n\geq n_0$, and, at second, $\lambda _n(S_n)\leq (1+{{\epsilon }\over 2})\lambda _{\cal F}(Q)$ for any $n\geq n_0$ according to the corollary 4.11. Therefore there is the set $S=S_{n_0}\times (B_{n_0+1}\times B_{n_0+2}\times \cdots)$ such that $S\supset Q$,  $S\in {\cal R}_{\cal E}$ and $\lambda _{\cal E}(S)<(1+{{\epsilon }\over 2})^2\lambda _{\cal F}(Q)$. Hence the estimation from above for $\bar \lambda _{\cal E}(Q)$   is obtained.

\medskip

{\bf II}. Estimation from below.
Let $Q=\{ x\in E:\ (x,f_k)\in [-{{d_k}\over 2},{{d_k}\over 2}]\}$ and $Q_n=\{ x\in E:\ (x,f_k)\in  [-{{d_k}\over 2},{{d_k}\over 2}] \ \forall \ k=1,\ldots,n;\ (x,f_i)=0\ \forall \ i>n\}$. Then for every $n\in \mathbb N$ the set $p_n={\mathbf P}_{E_n}(Q_n)$ is the orthogonal projection of the set  $Q_n$ onto the hyperplane  $E_n$.

Hence $p_n
=\Bigl\{ \sum\limits_{j=1}^n\sum\limits_{k=1}^n d_kt_kc_{j,k}e_j,\  t_1,\ldots ,t_n\in [-{{1}\over 2},{{1}\over 2}]\Bigr\}$ for any $n\in \mathbb N$. Then for any $x\in p_n$ and any $l=1,\ldots ,n$ the equality  $(x,f_l)=\sum\limits_{j=1}^n\sum\limits_{k=1}^n d_kt_kc_{j,k}c_{j,l}$ holds. Since $\sum\limits_{j=1}^{\infty } c_{j,k}c_{j,l}=\delta _{kl}$ then $(x,f_l)=d_lt_l- \sum\limits_{k=1}^n\sum\limits_{j=n+1}^{\infty } d_kt_kc_{j,k}c_{j,l},\ l=1,\ldots ,n$. 
Therefore 
\begin{equation}\label{maks} 
\sup\limits_{x\in p_n}|(x,f_l)-d_lt_l|\leq {1\over 2}\sum\limits_{k=1}^nd_k\sum\limits_{j=n+1}^{\infty } (|c_{j,k}||c_{j,l}|) \leq  {1\over 2}D_0\left ( \sum\limits_{k=1}^n\sum\limits_{j=n+1}^{\infty } (|c_{j,k}|) \right )\sum\limits_{j=n+1}^{\infty } |c_{j,l}|.
\end{equation}

Let us consider the numerical sequences $ \Delta _{l,n},\ n\in {\mathbb N},\, l\in \{ 1,...,n\}$ (see (\ref {delt})) and 
\begin{equation}\label{gam}
\gamma _n= \sum\limits_{k=1}^n\sum\limits_{j=n+1}^{\infty } |c_{j,k}| =\sum\limits_{k=1}^n\Delta _{k,n}.
\end{equation}

According to (\ref{maks}) 
\begin{equation}\label{maksi}
\sup\limits_{x\in p_n}|(x,f_l)|
\leq d_lt_l+{1\over 2}D_0\gamma _n\Delta _{l,n}\quad \forall \ l=1,...,n.
\end{equation}

Consider the set $\sigma _n\subset p_n$ where
\begin{equation}\label{base}
\sigma _n=\Bigl\{ \sum\limits_{k,j=1}^nd_kt_kc_{j,k}e_j,\; |t_k|\leq {1\over 2}\Bigl(1-{2\over d_0}\Delta _{k,n}\Bigr),\; k=1,\ldots,n \Bigr\}.
\end{equation} 
The set  $\sigma _n$ is the convex polyhedron in the space $E_n$. It is the image of $n$-dimensional rectangle $\Bigl\{ \sum\limits_{k,j=1}^nd_kt_kf_k;\ |t_k|\leq {1\over 2}(1-{2\over d_0}\Delta _{k,n}),\, k=1,\ldots,n \Bigr\}$ from the space  $F_n$ under the action of projector ${\bf P}_{F_n,E_n}$. Therefore 
\begin{equation}\label{LB}
\lambda _n(\sigma _n)=\det (C_n)\prod _{j=1}^n\Bigl( 1-{2\over d_0}\Delta _{j,n}\Bigr) \lambda _n(Q_n).
\end{equation}

The $n$-dimensional set $\sigma _n$ admits the extension $s_n$ along the directions $e_{n+1},e_{n+2},...$ such that $s_n\subset Q$.
Consider the set
\begin{equation}\label{vip}
s_n=\Bigl\{ \sigma _n+\sum\limits_{j=n+1}^{\infty }t_je_j ,\, t_j\in [-a_j,a_j],\, j\geq n+1\Bigr\}.
\end{equation}
Here for every $j\geq n+1$ the number $a_j$ is chosen from the interval  $(0,{1\over 2})$ such that the condition $s_n\subset Q$ holds. The condition $s_n\subset Q$  is equivalent to the system of inequalities  
$\sup\limits_{x\in s_n}|(x,f_k)|\leq  {{d_k}\over 2},\, \,k\in \mathbb N$. 

Since  $|t_l|\leq {1\over 2}(1-{2\over d_0}\Delta _{l,n}), $  for every $l=1,\ldots,n$ in the parametrization of the set (\ref{base}) then according to (\ref{maksi}) we obtain the estimates
 $\sup\limits_{x\in \sigma _n}|(x,f_l)|
\leq {1\over 2}d_l(1-\Delta _{l,n}({2\over {d_0}}-{{D_0}\over {d_0}}\gamma _n))$.

Since $a_j\leq {1\over 2} \ \forall \ j\geq n+1$ then $ \sup\limits_{x\in s_n}|(x,f_k)|\leq \sup\limits_{y\in \sigma_n}|(y,f_k)|+\sum\limits_{j=n+1}^{\infty }a_j|c_{jk}|\leq {1\over 2}d_k(1-\Delta _{k,n}({2\over {d_0}}-{{D_0}\over {d_0}}\gamma _n)) +{1\over 2}\Delta _{k,n}\leq {1\over 2}d_k(1-\Delta _{k,n}({1\over {d_0}}-{{D_0}\over {d_0}}\gamma _n))$
for any  $k\in \{ 1,\ldots,n\}$  according to (\ref{base}). Therefore the inequality $\sup\limits_{x\in s_n}|(x,f_k)|\leq {1\over 2}d_k,\, k=1,\ldots ,n$ holds for any sufficiently large $n$.

For every $k\geq 1+n$ we have  $\sup\limits_{x\in s_n}|(x,f_k)|\leq a_kc_{k,k}+{1\over 2}\sum\limits_{j\neq k}c_{jk}=a_k\alpha _{k}+{1\over 2}(l_k-\alpha _k)$. 
Hence conditions 
$\sup\limits_{x\in s_n}|(x,f_k)|< {1\over 2}d_k,\quad k>n+1$ are satisfied if the values $a_k$ are defined by the equalities  
$$
a_k={1\over {2\alpha _k}}(d_k-l_k+\alpha _k)={1\over 2}-{{l_k-d_k}\over {2\alpha _k}},\ k=n+1,...\ .
$$

Thus if the number $n\in \mathbb N$ is sufficiently large then there is the set $s_n\in {\cal R}_{\cal E}$ of type (\ref{vip}) such that  $Q\supset s_n$ and its $\lambda _{\cal E}$-measure has the estimates 
\begin{equation}\label{00}
\lambda _{\cal E}(s_n)=\lambda _n(\sigma _n)\prod _{j=n+1}^{\infty }2a_j\geq  [\det (C_n)\lambda _{n}(Q_n)\prod _{j=1}^n(1-{2\over d_0}\Delta _{l,n})
]\prod _{j=n+1}^{\infty }\Bigl(1-{{l_k-d_k}\over {\alpha _k}}\Bigr).
\end{equation}

According to corollary 4.6 we have  $\lim\limits_{n\to \infty }\det C_n=1$;  $\lim\limits_{n\to \infty }\lambda _{n}(Q_n)=\lambda _{\cal F}(Q)$ since $Q\in {\cal K}_{\cal F}$.

According to lemma 3.5 and theorem 4.2 we have $\lim\limits_{n\to \infty }(\sum\limits_{k=n+1}^{\infty }\sum\limits_{j=k+1}^{\infty }|c_{k,j}|+\sum\limits_{j=n+1}^{\infty }\sum\limits_{k=j+1}^{\infty }|c_{k,j}|)=0$ and  $\lim\limits_{n\to \infty }(\sum\limits_{l=1}^n\sum\limits_{j=n+1}^{\infty}
|c_{l,j}|)=0$.
Hence  $\lim\limits_{n\to \infty }(\sum\limits_{l=1}^n \Delta _{l,n})=0$.

Therefore, $\lim\limits_{n\to \infty }\left( \prod _{l=1}^n(1-{2\over {d_0}}\Delta _{l,n})\right) =
1$.

Since $\alpha _k\in (0,1]$ and $\lim\limits_{k\to \infty}\alpha _k =1$ then the series $\sum\limits_{k=1}^{\infty }{{l_k-1}\over {\alpha _k}}$ converges according to the condition (\ref{crit'}). According to the condition  $\lambda _{\cal F}(Q)\in (0,+\infty )$ the series $\sum\limits_{k=1}^{\infty }{{1-d_k}\over {\alpha _k}}$ converges.
Hence $\lim\limits_{n\to \infty }\prod\limits_{j=n+1}^{\infty }\Bigl(1-{{l_k-d_k}\over {\alpha _k}}\Bigr)=1$.
Therefore, $\lim\limits_{n\to \infty }\lambda _{\cal E}(s_n)=\lambda _{\cal F}(Q)$ according to (\ref{00}) and we obtain the estimate from below.
\hfill$\Box$

\bigskip

\section{The ring ${\cal R}_{\cal E}\bigcap {\cal R}_{\cal E}$ in the case of distant bases ${\cal F}$ and ${\cal E}$}.

{\bf Theorem 5.1}.
{\it If the condition (\ref{crit}) is not satisfied then for ane set $A\in {\cal R}_{\cal E}\bigcap {\cal R}_{\cal F}$ the equalities ${{\lambda _{\cal F}}}(A)=0$ and $\lambda _{\cal E}(A)=0$ hold.}

{\bf Lemma 5.2}. {\it Let $s_{\cal E}=\Pi ^0\backslash (\bigcup\limits_{i=1}^n\Pi ^i),\
\Pi ^k\in {\cal K}_{\cal E}\ \forall \ k=0,\ldots, n$, и $\lambda _{\cal E}(s_{\cal E})>0$. 
Then there is the rectangle $P _{\cal E}\in {\cal K}_{\cal E}$ such that $P _{\cal E}\subset s_{\cal E}$ and $\lambda _{\cal E}(P_{\cal E})>0$.}

{\it Proof}. Firstly we consider the case 
$s_{\cal E}=\Pi \backslash \Pi '$ where $\Pi ,\Pi '\in {\cal K}_{\cal E}$ (we can assume that $\Pi '\subset \Pi$). Let $[a_j,b_j),\, [a_j',b_j'),\ j\in {\mathbb N}$ be the edges of the rectangles $\Pi ,\Pi '$ that are collinear to the vector $e_j$ for every $ j\in {\mathbb N}$. Then $[a_j,b_j)\supset [a_j',b_j')$. Since $\lambda _{\cal E}(\Pi \backslash \Pi ')>0$ then there in the number $k\in {\mathbb N}$ such that $[a_k,b_k)\backslash [a_k',b_k')\supset [a_k'',b_k'')$ and $b_k''-a_k''>0$. Let $\Pi ''$ be the rectangle such that $k$-th edge of $\Pi''$ is the segment $[a_k'',b_k'')$ and its $j$-th edge coincides with $j$-th edges of the rectangle $\Pi$ for every $j\neq k$. Then the rectangle $\Pi ''$ satisfies the following conditions: $\Pi ''\subset {\cal K}_{\cal E}$, $\Pi ''\subset s_{\cal E}$ and $\lambda _{\cal E}(\Pi '')>0$. 

In general case $s_{\cal E}=\Pi ^0\backslash (\bigcup\limits_{i=1}^n\Pi ^i),\ \Pi ^i\in {\cal K}_{\cal E},\, i=0,1,...,n$, the statement of the lemma 5.2 can be obtained by the applying of the induction method with respect to $n\in \mathbb N$.
$\hfill \Box$

\medskip

Let $ l_j=\sum\limits_{k=1}^{\infty }|c_{k,j}|$ be the lengths of the projection of the unit rectangle $Q\in {\cal K}_{\cal F}$ onto the axis $Oe_j$ for every $j\in \mathbb N$. Let  $l_k^T=\sum\limits_{j=1}^{\infty }|c_{k,j}|$  be the lengths of the projection of the unit rectangle $\Pi \in {\cal K}_{\cal E}$ onto the axis $Of_k$ for every $k\in \mathbb N$. Then $l_j\geq 1,\ l^T_k\geq 1$ according to lemma 3.2. Since the condition (\ref{crit}) is not satisfied then the following two conditions hold 
$$
\prod\limits_{j=1}^{\infty } l_j=+\infty ;\quad \prod\limits_{k=1}^{\infty } l^T_k=+\infty .
$$ 
These two condition are equivalent to the following two equalities \begin{equation}\label{ryad}
\sum\limits_{j=1}^{\infty }(l_j-1)=+\infty ;\quad \sum\limits_{k=1}^{\infty }(l^T_k-1)=+\infty .
\end{equation}
Let $L^T=\sup\limits_{k\in {\mathbb N}}l_k^T$ and $L=\sup\limits_{j\in {\mathbb N}} l_j$. Then there are three possible cases for the sequences 
$\{ l^T_k\}$ and $\{  l_j\}$.

1)  $L^T=+\infty $ и $L=+\infty $. 

2)  $L^T<+\infty $ и $L<+\infty $. 

3) either $L^T=+\infty $ and $L<+\infty $, or vice versa.
 
Let us prove the statement of the theorem 5.1 for any of these three cases. 

\medskip

1) Consider the case $L^T=+\infty $ and $L=+\infty $. 

{\bf Lemma 5.3.} {\it Let $L^T=+\infty $. If $\Pi \in {\cal K}_{\cal E}$ and $\lambda _{\cal E}(\Pi)>0$ then  ${\overline {\lambda _{\cal F}} }(\Pi)=+\infty $. 

Let  $L=+\infty $. If $Q \in {\cal K}_{\cal F}$ and $\lambda _{\cal F}(Q)>0$ then ${\overline {\lambda _{\cal E}} }(Q)=+\infty $.}

{\bf Proof}.
Since $L^T=+\infty $ then there are two possible cases:   
either
{\bf i)} $\exists \ k_0\in {\mathbb N}:\  l^T_{k_0}=+\infty$;
or 
{\bf ii)} $\forall \ k\in {\mathbb N}:\ l^T_k<+\infty $.

Let us prove that in every of these cases the condition $\Pi \in {\cal K}_{\cal E}:\ \lambda _{\cal E}(\Pi )>0$ implies   ${\overline {\lambda _{\cal F}} }(\Pi )=+\infty $.

{\bf i)}. In the case i) it is impossible to cover a rectangle  $\Pi \in {\cal K}_{\cal E}:\ \lambda _{\cal E}(\Pi )>0$ by the finite union $S$ of rectangles from the collection ${\cal K}_{\cal F}$.  In fact, the projection of the rectangle $\Pi$ onto the the axis $Of_{k_0}$ is unbounded segment since $l_{k_0}^T=+\infty $. However, the projection of a set $S$ onto the axis $Of_{k_0}$ is the finite union of bounded segments. Therefore, the statement is proved and ${\overline {\lambda _{\cal F}} }(\Pi )=+\infty $.

{\bf ii)}. 
If $\Pi \in {\cal K}_{\cal E}:\ \lambda _{\cal E}(\Pi )>0$ then the sequence of lengths of projections of the rectangle $\Pi$ onto the axes $Of_k,\, k\in {\mathbb N}$ is unbounded. Because if $\{r_k\}$ is the sequence of lengths of edges of the rectangle $\Pi$ then $\inf\limits_jr_j=r_0>0$ according to the condition $\lambda _{\cal E}(\Pi)>0$. Therefore the length of the projection of the rectangle $\Pi $ onto the axis $Of_k$ no less than $r_0 l^T_k$. 
The sequence $l_k^T$ is unbounded according to {\bf ii).} Hence  the sequence of lengths of projections of the rectangle $\Pi$ onto the axes $Of_k,\, k\in {\mathbb N}$ is unbounded.

On the other hand, if $S=\bigcup\limits_{k=1}^mQ_k$ is the finite union of rectangles from the collection ${\cal K}_{\cal F}$ then the sequence of projections of the set $S$ onto an axis $Of_j$ is bounded set. Hence in the case {\bf ii)} it is impossible to cover the rectangle $\Pi \in {\cal K}_{\cal E}:\ \lambda _{\cal E}(\Pi )>0$ by the finite union $S$ of rectangles from the collection ${\cal K}_{\cal F}$. Therefore, ${\overline {\lambda _{\cal F}} }(\Pi )=+\infty $. $\hfill \Box$

\medskip

{\bf Corollary 5.4}. {\it If  $L=+\infty $ then the the condition $A\in {\cal R}_{\cal E}\bigcap {\cal R}_{\cal F}$ implies  $\lambda _{\cal F}(A)=0$. If $L^T=+\infty $ then the condition $A\in {\cal R}_{\cal E}\bigcap {\cal R}_{\cal F}$ implies $\lambda _{\cal E}(A)=0$. If the condition 1) is satisfied then for any set  $A\in {\cal R}_{\cal E}\bigcap {\cal R}_{\cal F}$ the equalities  $\lambda _{\cal F}(A)=\lambda _{\cal E}(A)=0$ hold.}

{\bf Proof}.
Let $A\in {\cal R}_{\cal E}\bigcap {\cal R}_{\cal F}$. 
Let's assume the opposite that $\lambda _{\cal E}(A)>0$. Then according to lemma 5.2 there is a rectangle $p\in {\cal K}_{\cal E}$ such that $p\subset A$ and $\lambda _{\cal E}(p)>0$. Then according to the lemma 5.3 the condition  $L^T=+\infty$ implies that ${\overline {\mu _{\cal F}} }(A)\geq {\overline {\mu _{\cal F}} }(p)=+\infty $. This is the contradiction to the condition $A\in {\cal R}_{\cal F}$. Therefore,  $\lambda _{\cal E}(A)=0$. Analogously, the condition $ L=+\infty$ implies that $\lambda _{\cal F}(A)=0$. $\hfill \Box$

\medskip

2) Consider the case $L^T<+\infty $ and $L<+\infty $. 

{\bf Lemma 5.5}. {\it Let the condition 2) be satisfied. Let $\Pi \supset Q$ where $\Pi \in {\cal K}_{\cal E},\, Q\in {\cal K}_{\cal F}$. Then $\lambda _{\cal F}(Q)=0$.}

{\bf Proof}. Let's assume the opposite that L  $\lambda _{\cal F}(Q)>0$. 
Let $\{ d_k\}$ be the sequence of lengths of the edges of the rectangle $Q$. We can assume that  $d_k\leq 1$. In fact, in opposite case we can change the inscribed onto the smaller inscribed rectangle $Q':\ Q'\subset Q\subset \Pi$ such that  lengths of the edges of the rectangle $Q'$ no greater than $1$. Then $\delta _k=1-d_k\geq 0 \ \forall \ k$, 
$\lim\limits_{k\to \infty }d_k=1$ and 
\begin{equation}\label{l-fin}
\sum\limits_{k=1}^{\infty }\delta _k<+\infty 
\end{equation} 
according to the condition $\lambda _{\cal F}(Q)\in (0,+\infty )$.

Therefore, if $\Pi \in {\cal K}_{\cal E}$, $\Pi \supset Q$ and $\{ D_j\}$ is the sequence of lengths of edges of the rectangle $\Pi$ then $$D_j\geq \sum\limits_{k=1}^{\infty }d_k|(f_{k},e_{j})|=
\sum\limits_{k=1}^{\infty }(1-\delta _k)C_{jk}= l_j-\sum\limits_{k=1}^{\infty }\delta _kC_{jk}.$$ Hence 
$$
\sum\limits_{j=1}^{\infty }(D_j-1)\geq 
\sum\limits_{j=1}^{\infty }( l_j-1 - \sum\limits_{k=1}^{\infty }C_{jk}\delta _k)= \sum\limits_{j=1}^{\infty }( l_j-1)-\sum\limits_{k=1}^{\infty }(\sum\limits_{j=1}^{\infty }C_{jk})\delta _k=\sum\limits_{j=1}^{\infty }( l_j-1)-\sum\limits_{k=1}^{\infty } l^T_k\delta _k\geq $$
$$
\geq \sum\limits_{j=1}^{\infty }(l_j-1)- L^T\sum\limits_{k=1}^{\infty }\delta _k=+\infty
$$ 
since the first series diverges according to the assumption of violation of the condition  (\ref{crit}) and the second series is converges according to the condition (\ref{l-fin}) with finiteness of the value $L^T=\sup\limits_k l^T _k$. 
But the condition $\Pi \in {\cal K}_{\cal E}$ implies that $\sum\limits_{j=1}^{\infty }(D_j-1)<+\infty$. The obtained contradiction proves the statement of lemma 5.5. $\hfill \Box$

{\bf Corollary 5.6.} {\it Let the condition 2) is satisfied. Then for any rectangle $\Pi \in {\cal K}_{\cal E}$ the equality ${\underline { \lambda _{\cal F}}}(\Pi )=0$ holds.}

{\bf Proof.} Let's assume the opposite that  there is a rectangle $\Pi \in {\cal K}_{\cal E}$ such that  ${\underline { \mu _{\cal F}}}(\Pi )>0$. Then there is the set $s\in r_{\cal F}$  such that $\lambda _{\cal F}(s)>0$ and $s\subset \Pi$. Hence according to lemma 5.2 there is a rectangle  $q\in {\cal K}_{\cal F}$ such that $q\subset s\subset \Pi$ and $\lambda _{\cal F}(q)>0$. Therefore $\lambda _{\cal E}(\Pi)=+\infty $ according to lemma 5.5. This is the contradiction with the condition $\Pi \in {\cal K}_{\cal E}$. $\hfill \Box$

\medskip

Let us prove that it is impossible to cover the rectangle $q\in {\cal K}_{\cal F}$ with positive measure $\lambda _{\cal F}(q)>0$ by the finite union $S\in r_{\cal E}$ of rectangles $\Pi _1,...,\Pi _m\in {\cal K}_{\cal E}$. Then we show that $\lambda _{\cal F}(A)=0$ for any  $A\in {\cal R}_{\cal E}\bigcap {\cal R}_{\cal F}$.

{\bf Lemma 5.7}. {\it 
Let $S=\bigcup\limits_{i=1}^m\Pi _i$ where $\Pi _i\in {\cal K}_{\cal E},\, i=1,...,m$ and $\Pi _i\neq \varnothing \ \forall \ i=1,...,m$. Then there is a hyperplane $\Gamma$ of finite codimension $k\leq m-1$ of type  
\begin{equation}\label{hypergama}
\Gamma _{c_1,\ldots,c_k}^{j_1,\ldots,j_k}=\{ x\in E:\ (x,e_j)=c_j,\, j=1,\ldots , k\}
\end{equation} 
such that $\Gamma \bigcap S=\Gamma \bigcap \Pi _{j_*}$ for some $j_*\in \{1,\ldots , m\}$.}

{\bf Proof.} Let us prove the statement by the induction. Firstly we note that the statement is true for $m=1$ since in this case $S=\Pi _1$ and $\Gamma =E$. 

Let ${\bf P}_{e_j}^{\bot }$ be the orthogonal projection onto the subspace $(x,e_j)=0$ in the space $E$ for every $j\in \mathbb N$. For every $j\in \mathbb N$ the intersection $\Gamma _{j,c}\bigcap \Pi _1$ is either empty set or rectangle in the hyperplane $\Gamma _{j,c}$. 

Let $m\in \mathbb N$. 

The intersection $\Gamma _{j,c}\bigcap S=\bigcup\limits_{i=1}^m(\Gamma _{j,c}\bigcap \Pi _i)$ is the union of no more than $m$ nonempty rectangles $\Gamma _{j,c}\bigcap \Pi _i,\, i=1,...,m$ in the hyperplane $\Gamma _{j,c}$.  
Therefore for every $j\in \mathbb N$ the set-valued function ${\bf P}^{\bot }_{e_j}(\Gamma _{j,c}\bigcap S),\ c\in \mathbb R$, can has only finite number of 
values and these values lies in a set of subsets of the subspace $E\ominus {\rm span}(e_j)$.

If for any $j\in \mathbb N$ the set-valued function ${\bf P}^{\bot }_{e_j}(\Gamma _{j,c}\bigcap S),\ c\in \mathbb R$, has only one nonempty value $B_j$, then the set $S$ is the rectangle $(\times) _{j=1}^{\infty }{\bf P}_{Oe_j}(S)$
and the statement of the lemma is true for hyperplane  $\Gamma _{j,c}$ for any $j$ and for any $c\in {\bf P}_{Oe_j}(S)$. 

In opposite case there is a number  $j\in \mathbb N$
and there are $c',c''\in {\bf P}_{Oe_j}(S)$
such that the sets
$\Gamma _{j,c'}\bigcap S$ and $\Gamma _{j,c''}\bigcap S$ are nonempty sets and ${\bf P}^{\bot }_{e_j}(\Gamma _{j,c'}\bigcap S)\neq {\bf P}^{\bot }_{e_j}(\Gamma _{j,c''}\bigcap S)$. 
Hence 
for at least one of two numbers $c',c''$ (for certainty, for the number $c'$) the set ${\bf P}^{\bot }_{e_j}(\Gamma _{j,c'}\bigcap S)$ can't be the union of $m$ nonempty rectangles  ${\bf P}^{\bot }_{e_j}(\Gamma _{j,c'}\bigcap \Pi _i),\, i=1,...,m$ in the space $E\ominus {\rm span}(e_j)$. Therefore there is the numbers $j\in \mathbb N$ and $c'\in \mathbb R$ such that the set $\Gamma _{j,c'}$ is the union of no more than $m-1$ rectangles in the space $E\ominus {\rm span}(e_j)$.  Then we can apply the assumption of induction to the set $\Gamma _{j,c'}(S)$  in the space $E\ominus {\rm span}(e_j)$. Thus we obtain the statement of the lemma. $\hfill \Box$

{\bf Corollary 5.8.} {\it
Let $q\in {\cal K}_{\cal F}$ and $q\subset S$ where $S=\bigcup\limits_{i=1}^m\Pi _i,\ \Pi _i\in {\cal K}_{\cal E} \ \forall \ i=1,...,m.$ 
Then there are a hyperplane $\Gamma $ of type (\ref{hypergama}) and a number $i^*\in \{ 1,...,m\}$ such that $\Gamma \bigcap q\subset \Gamma \bigcap \Pi _{i_*}$.
}

{\bf Lemma 5.9}. {\it Let ${\cal E},\, {\cal F}$ be a pair of ONB such that the condition (\ref{crit}) is not satisfied. If $q\in {\cal K}_{\cal F}$ and $\lambda _{\cal F}(q)>0$ then ${\overline {\lambda _{\cal E}}}(q)=+\infty $. If $p\in {\cal K}_{\cal E}$ and $\lambda _{\cal E}(p)>0$ then ${\overline {\lambda _{\cal F}}}(p)=+\infty $}.

{\bf Proof.} We can assume without loss of generality that $q$ is a centered rectangle such that the sequence of length of its edges $\{ d_k\}$ satisfies the condition $d_k\leq 1 \, \forall \, k$ (see proof of the lemma 5.5).
Since $\lambda _{\cal F}(q)>0$ then  $d_0=\inf\limits_{k}d_{k}>0$ and $D_0=\sup\limits_{k}d_k<+\infty $.  The lengths of the projections of the rectangle $q$ onto axes 
$Oe_j$ form the numerical sequence $\{ D_j\}$. Then $D_j\leq D_0 L$. On the other hand, the condition 
\begin{equation}\label{inf-5}
\prod\limits_{j=1}^{\infty }D_j=\infty 
\end{equation} 
holds according to lemma 5.5.
For every $j\in \mathbb N$ the projection of the rectangle $q$ onto the axis $Oe_j$ is the segment containing the interval $\Bigl(-{{D_j}\over 2}, {{D_j}\over 2}\Bigr)$. 

The rectangle $q$ admits the parametrization $$q=\Bigl\{ x=\sum\limits_{i=1}^{\infty }d_it_if_i,\ t_i\in [-{1\over 2},{1\over 2}], \, \sum\limits_{i=1}^{\infty }t_i^2<+\infty \Bigr\}.
$$ 

Let us fix a numbers $j\in \mathbb N$ and $b \in (-{{D_j}\over 2},{{D_j}\over 2})$. Consider the intersection  $\Gamma _{j,b}\bigcap q= \{ x\in q:\ x_j=b\}$. The set $\Gamma _{j,b}\bigcap q$ admits the parametrization 
\begin{equation}\label{brus}
\Gamma _{j,b}\bigcap q=\Bigl\{ x=\sum\limits_{i=1}^{\infty }d_it_if_i,\ t_i\in \Bigl[-{1\over 2},{1\over 2}\Bigr], \, \sum\limits_{i=1}^{\infty }t_i^2<+\infty ;\ \sum\limits_{i=1}^{\infty }d_ic_{i,j}t_i=b \Bigr\} .
\end{equation}
Since $|b|<{{D_j}\over 2}$ then $\Gamma _{j,b}\bigcap q\neq \emptyset$. Hence there is a vector $\{ t_i^0\}\in \ell _2$ such that  $\sum\limits_{i=1}^{\infty }d_ic_{i,j}t_i^0=b$ and $|t_i^0|<{1\over 2}\ \forall \ i\in \mathbb N$.

Since  $\alpha _j=\max\limits_{i\in {\mathbb N}}|c_{i,j}|>0$ 
for given $j$ then we can define the positive number 
$$
\delta _j=\min \{ {1\over {2\alpha _j d_j}}({{D_j}\over 2}-|b| ),\, {1\over 2}-|t_j^0|\}>0.
$$ 
Hence for given number $j$ there is the number $i_0>j$ such that $\sum\limits_{i=i_0}^{\infty }d_i|c_{i,j}|<{1\over {2}}\alpha _jd_j\delta _j
$.
Therefore for every collection of numbers $\hat t_j=(t_1,...,t_{j-1},t_{j+1},...)$ such that
$$
t_i\in \Bigl[-{{d_j\delta _j}\over {4i_0D_0}},{{d_j\delta _j}\over {4i_0D_0}}\Bigr] ,\ i=1,\ldots,i_0,\, i\neq j;\quad t_i\in \Bigl[-{1\over 2},{1\over 2}\Bigr],\ i>i_0
$$ 
there is the number $t_j=t_j(\hat t_j)\in (t_j^0-\delta _j, t_j^0+\delta _j)$ such that $d_jc_{j,j}t_j+\sum\limits_{i\neq j}d_ic_{i,j}t_i=b$. 

Thus,  the points with parameters $t_i,\, i\in \mathbb N$, in the parametrization  (\ref{brus}) belong to the set $\Gamma _{j,b}\bigcap q$ if the following conditions hold:
$$
t_i\in [-{1\over 2},{1\over 2}],\ i> i_0;$$
\begin{equation}\label{t-sect}
|t_i-t_i^0|< {{d_j\delta _j}\over {4i_0D_0}},\, i=1,\ldots,i_0,\, i\neq j;
\end{equation}
$$
t_j=t_j(\hat t_j)\in (t_j^0-\delta _j, t_j^0+\delta _j.)
$$ 
For any $x\in \Gamma _{j,b}\bigcap q$ which is given by parametrization (\ref{brus}) with the parameters from the set (\ref{t-sect})
we have $(x,e_k)=
\sum\limits_{i=1}^{\infty }d_ic_{i,k}t_i$. Hence
$$\sup\limits_{x\in \Gamma _{j,b}\bigcap q}|(x,e_k)|\geq \sum\limits_{i=1}^{\infty}{1\over 2}d_i|c_{i,k}|-\sum\limits_{i=1}^{i_0}{1\over 2}d_i|c_{i,k}|$$ for any $k\in \mathbb N$. Therefore,
for every $k>i_0$ we obtain 
$$
\{ (x,e_k),\ x\in \Gamma _{j,b}\bigcap q\} \supset \Bigl(-{1\over 2}D_k+{1\over 2}\sum\limits_{i=1}^{i_0}d_i|c_{ik}|,{1\over 2}D_k-{1\over 2}\sum\limits_{i=1}^{i_0}d_i|c_{ik}|\Bigr).
$$ 

Since $\sum\limits_{k=k_0+1}^{\infty }(D_k-1)-\sum\limits_{k=k_0+1}^{\infty }\sum\limits_{i=1}^{i_0}d_i|c_{ik}|\geq \sum\limits_{k=k_0+1}^{\infty }(D_k-1) -D_0\sum\limits_{i=1}^{i_0} l^T_i =+\infty $  according to (\ref{inf-5}) and the condition $L^T<\infty $. Then  we obtain 
\begin{equation}\label{produ}
\prod\limits _{k=k_0+1}^{\infty }(D_k-\sum\limits_{i=1}^{i_0}d_i|c_{ik}|)=+\infty .
\end{equation}

Therefore every intersection of the rectangle $q$ by the one-codimensional hyperplane of type  $\Gamma _{j,c}=\{ x:\ (x,e_j)=c,\, j\in {\mathbb N},\, c\in (-{{D_j}\over 2},{{D_j}\over 2})\}$ can not be covered by one measurable rectangle 
$\Pi ^*\in {\cal K}_{\cal E}$. Because if $\Pi ^* \supset (\Gamma _{j,b}\bigcap q)$ then for every $k>i_0$ the $k$-th edge of rectangle $\Pi ^*$ should has the length no less than
$D_k-\sum\limits_{i=1}^{i_0}d_i|c_{ik}|$. Hence the condition $\Pi ^*\in {\cal K}_{\cal E}$ is violated according to (\ref{produ}).

The same reasoning allows us to show that
every intersection of rectangle $q$ by the $m$-codimensional hyperplane $
\Gamma ^{j_1,\ldots,j_m}_{c_1,\ldots ,c_m}$ of type  
$$
\{ x:\ (x,e_{j_1})=c_1,\ldots, (x,e_{j_m})=c_m,\, j_1,\ldots,j_m\in {\mathbb N},\, c_k\in (-{{D_k}\over 2},{{D_k}\over 2})\ \forall \ k\in \{ j_1,\ldots,j_m\}\}$$  with some $m\in \mathbb N$ can not be covered by one measurable rectangle $\Pi ^*\in {\cal K}_{\cal E}$. 

Let's assume that $\overline \lambda _{\cal E}(q)<+\infty$. Then there is the set $S=\bigcup\limits_{k=1}^N\Pi _k,\ N\in {\mathbb N},\, \Pi _k\in {\cal K}_{\cal E}\, \forall \, k\in \{ 1,...,N\}$, such that $q\subset S$. Then ${x\in S:\ (x,e_j)=c }\neq \emptyset $ for every $j\in {\mathbb N}$ and for every $c\in (-{1\over 2}D_j,{1\over 2}D_j)$.
Therefore,
according to the corollary 
5.8 the condition  $q\subset S=\bigcup\limits_{k=1}^N\Pi _k$ implies that there are numbers $j^*\in \{1,...,N\}$ and $c_1,...,c_N$ such that  $c_i\in (-{1\over 2}D_i,{1\over 2}D_i)$ and 
$\Gamma ^{j_1,\ldots, j_m}_{c_1,\ldots ,c_m}\bigcap q\subset \Pi _{j^*}$. Hence the intersection $q\bigcap \Gamma ^{j_1,\ldots, j_m}_{c_1,\ldots ,c_m}$ is covered by one rectangle $\Pi _{j^*}$.

The obtained contradiction prove that the rectangle $q\in {\cal K}_{\cal F}:\ \lambda _{\cal F}(q)>0$ can't be covered by the finite union of rectangles from the collection ${\cal K}_{\cal E}$. Therefore,  ${\overline {\lambda _{\cal E}}}(q)=+\infty $. $\hfill \Box$

{\bf Corollary 5.10.}  {\it Let the pair of bases satisfies the condition 2). Then for any $A\in {\cal R}_{\cal E}\bigcap {\cal R}_{\cal F}$ the equalities $\lambda _{\cal E} (A)=\lambda _{\cal E} (A)=0$ hold.} (The proof is the same as the proof for the corollary 5.4).

\medskip

3) Consider the case $L^T=+\infty $, $L<+\infty $, or vice versa. 

Let $L^T=+\infty $, $L<+\infty $.
Then the sequence $l^T_k,\ k\in \mathbb N$, either is unbounded real valued sequence or takes values $+\infty $. The sequence $ l_j,\ j\in \mathbb N$ is bounded against, but $\sum\limits_{j=1}^{\infty }(l_j-1)=+\infty $. 
Since $L^T=+\infty$ then according to the corollary 5.4 $\lambda _{\cal E}(A)=0$ for any set $A\in {\cal R}_{\cal E}\bigcap {\cal R}_{\cal F}$.

{\bf Theorem 5.11}. {\it Let the condition 3) is satisfied. Then $\lambda _{\cal F}(A)=0$ for any set $A\in  {\cal R}_{\cal E}\bigcap {\cal R}_{\cal F}$}.

{\bf Proof}. Let's assume the opposite that  there is a set  $A\in  {\cal R}_{\cal E}\bigcap {\cal R}_{\cal F}$ such that $\lambda _{\cal F}(A)>0$. Then according to the lemma 5.2 there is the rectangle $q\in {\cal K}_{\cal F}$ such that  $q\subset A$ and $\lambda _{\cal F}(q)>0$.
Let $\{ d_k\}$ be the sequence of lengths of edges of the rectangle $q\in {\cal K}_{\cal F}$. Then we can assume without loss of generality that   $d_k\leq 1 \, \forall \, k\in \mathbb N$. Since  $\lambda _{\cal F}(q)>0$  then $\sum\limits_{k=1}^{\infty }|d_k-1|<+\infty$. In particular 
\begin{equation}\label{limd}
\lim\limits_{k\to \infty }d_k=1.
\end{equation} 
The lengths of projections of the rectangle $q$ onto the axes $Oe_j,\, j\in \mathbb N$ are  $$D_j=\sum\limits_{k=1}^{\infty }|c_{j,k}|d_k,\  j\in \mathbb N.$$
Therefore, if $\Pi \in {\cal K}_{\cal E}$ and $\Pi \supset q$ then 
$$\lambda _{\cal E}(\Pi )\geq \exp (\sum\limits_{j=1}^{\infty }\ln (D_j)).$$

To prove the theorem 5.11 we firstly obtain the statements 5.12-5.16.

\medskip

{\bf Lemma 5.12}. {\it If the condition
\begin{equation}\label{upl}
{\overline {\lim\limits_{j\to \infty }}} l_j=L_0>1,
\end{equation} 
is satisfied then $\sum\limits_{n=1}^{\infty}\max \{ 0,\ln(D_j)\}=+\infty$.}

{\bf Proof}. Let's fix $\epsilon >0$. According to (\ref{limd})
there is the number $K_{\epsilon }$ such that  $d_k>1-\epsilon $ for any $k>K_{\epsilon }$.

Since $\lim\limits_{j\to \infty }c_{k,j}=0$ for any  $k$ then $\lim\limits_{j\to \infty }\left( \max\limits_{k\in 1,\ldots,K_{\epsilon }}c_{k,j}\right) =0$.
Hence there is the number $J_{\epsilon}$ such that $\sum\limits_{k=1}^{K_{\epsilon }}c_{k,j}<\epsilon $ for every $j> J_{\epsilon }$.
Therefore, for every $j>J_{\epsilon}$ we have the estimate  
\begin{equation}\label{abc}
D_j=\sum\limits_{k=1}^{\infty }c_{k,j}d_k\geq \sum\limits_{k=K_{\epsilon}}^{\infty }c_{k,j}(1-\epsilon )\geq (1-\epsilon )( l_j-\epsilon ).
\end{equation}

Since  
$\epsilon >0$ is arbitrary 
then according to (\ref{upl}) we can choose the value $\epsilon >0$ in (\ref{abc}) such that there is a strictly monotone sequence of numbers $\{ j_n\}$ satisfying following condition: $D_{j_n}>{1\over 2}(1+L_0)>1$ for any $n\in \mathbb N$. Hence $\sum\limits_{n=1}^{\infty }(D_{j_n}-1)=+\infty $ and the statement of lemma 5.12 is proved. $\hfill \Box$

{\bf  Corollary 5.13}. {\it If the condition (\ref{upl}) is satisfied then the rectangle  $q\in {\cal K}_{\cal F}$ with positive measure $\lambda _{\cal F}(q)>0$ can't be covered by the rectangle  
$\Pi \in 
{\cal K}_{\cal E}$.}

\medskip

Let us introduce the notation  $\delta _k=1-d_k,\, k\in \mathbb N$. 

{\bf Lemma 5.14}. {\it Let the condition (\ref{upl}) is violated. If $\| \delta \|_2={\sqrt {\delta _1^2+\delta _2^2+\ldots}}<1$ then the rectangle $q$ with positive measure $\lambda _{\cal F}(q)>0$ can't be inscribe into any rectangle $\Pi \in {\cal K}_{\cal E}$.
}

{\bf Proof}. Since the condition (\ref{upl}) is not satisfied then the limit $\lim\limits_{j\to \infty}l_j=1$ exists since $l_j\geq 1 \ \ \forall \ j\in \mathbb N$. 
According to lemma 3.2  the condition $ \alpha _j<{1\over {\sqrt 2}}$ implies the estimate $ l_j\geq {\sqrt 2}$;\ the condition  $\alpha _j>{1\over {\sqrt 2}}$ implies the inequality $l_j\geq \alpha _j+{\sqrt {1-\alpha _j^2}}$ (remember that $\alpha _j=\max\limits_{k\in {\mathbb N}}|c_{k,j}|,\ j\in \mathbb N$.
Hence the equality $\lim\limits_{j\to \infty} l_j=1$ implies that $\lim\limits_{j\to \infty}\alpha _j=1$. If $\alpha _j=1-\beta _j$ (see the notation used in the lemma 3.5) then $l_j\geq 1-\beta +{\sqrt {2\beta _j-\beta _j^2}}$. Therefore, the asymptotic equality 
\begin{equation}\label{omal}
\beta _j={1\over 2}{{(l_j-1)}^2}(1+o(1))
\end{equation} 
as $j\to \infty $ holds.
For every $j\in \mathbb N$ we have the estimate
$$
D_j-1=\sum\limits_{k=1}^{\infty }c_{k,j}d_k-1=l_j-1-\sum\limits_{k=1}^{\infty }c_{k,j}\delta _k=
$$
\begin{equation}\label{32}
=
l_j-1-\sum\limits_{k\neq j}^{\infty }c_{k,j}\delta _k-c_{j,j}\delta _j\geq 
l_j-1-\delta _j-{\sqrt {1-\alpha _j^2}}\| \delta \|_2.
\end{equation}

According to (\ref{omal}) we have ${\sqrt {1-\alpha _j^2}}\| \delta \|_2={\sqrt {\beta _j(1+\beta _j)}}\| \delta \|_2\sim (l_j-1)\|\delta _j\|_2$ as $j\to \infty $. Therefore, there are the numbers $J\in \mathbb N$ and $\sigma \in (0, 1-\| \delta \|_2)$ such that ${\sqrt {1-\alpha _j^2}}\| \delta \|_2\leq (1-\sigma )(l_j-1)$ for every $j>J$. Hence 
$$
D_j-1\geq \sigma (l_j-1)-\delta _j \ \ \forall \ \ j>J
$$
according to ({\ref {32}}).
Since the condition (\ref{crit}) is not satisfied then the nonnegative series  $\sum\limits_{j=1}^{\infty }(l_j-1)$ diverges.
Therefore
$
\sum\limits_{j=1}^{\infty }
(D_j-1)=+\infty ,
$ 
hence $\lambda (\Pi )=\exp (\sum\limits_{j=1}^{\infty } \ln(D_j))=+\infty $. $\hfill \Box$

{\bf Lemma 5.15}. {\it Let $q\in {\cal K}_{\cal F}$ and the sequence $\{ d_k\}$ of lengths of its edges satisfy the condition of the lemma  5.14. Then ${\overline {\lambda _{\cal E}}}(q)=+\infty$.}

The proof of the lemma 5.15 repeats the proof of the lemma 5.9.

{\bf Lemma 5.16}. {\it Let $q\in {\cal F}$ and $\lambda _{\cal F}(q)>0$. Then   ${\overline {\lambda _{\cal E}}}(q)=+\infty$.}

{\bf Proof.} Let's assume the opposite that  ${\overline {\lambda _{\cal E}}}(q)<+\infty$. Hence there is a set $S=\bigcup\limits_{s=1}^m\Pi _s$ such that  $\Pi _1,\ldots,\Pi _m\in {\cal K}_{\cal E}$ and $q\subset S$. 

Since $\lambda _{\cal F}(q)>0$ then there is a number $N\in \mathbb N$ such that $\sum\limits_{k=N+1}^{\infty }\delta _k<1$.

We can assume that the rectangle $q$ is centered.
Consider the orthogonal projections  $q_N={\bf P}_{F_N}(q),\ q^N={\bf P}_{F^N}(q)$ of the rectangle $q$ onto $N$-dimensional subspace $F_N={\rm span}(f_1,\ldots ,f_N)$ and onto its orthogonal complement respectively. Let $Q_N^1$ be the centered unit rectangle in the subspace $F_N$. Let $Q=Q_N^1\times q^N$. Then the rectangle $Q$
satisfies the conditions of lemma 5.14. Therefore, \begin{equation}\label{QQ}
{\overline {\lambda _{\cal E} }}(Q)=+\infty.
\end{equation}

Lengths of edges of the rectangle $q$ satisfy  conditions $d_k\leq 1 \, \forall \, k\in \mathbb N$. 
Hence $q_N\subset Q_N^1$ and $q\subset Q$. Since $\lambda _N(q_N)>0$ then there is the collection of vectors $h_1,\ldots,h_M\in F_N$ such that $\bigcup\limits_{j=1}^M{\bf S}_{h_j}(q_N)\supset Q_N^1$ (here ${\bf S}_{h_j}(q_N)=q_N+h_j$). Since $\bigcup\limits_{j=1}^M{\bf S}_{h_j}(q)\supset Q$ therefore $\bigcup\limits_{j=1}^M{\bf S}_{h_j}(S)\supset Q$. Thus we obtain the contradiction with the condition (\ref{QQ}).~$\hfill \Box$

Therefore, the statement of the theorem 5.11 is the consequence of the lemma 5.16.  $\hfill \Box$

{\bf Corollary 5.17}. {\it Let the condition 3) be satisfied. Then the equalities $\lambda _{\cal E} (A)=\lambda _{\cal E} (A)=0$ hold for any set $A\in {\cal R}_{\cal E}\bigcap {\cal R}_{\cal F}$.}

The theorem 5.1 follows from the corollaries 5.4, 5.10, 5.17.

\section{Translation and rotation-invariant measure}

Let $\mathcal S$ be a set of ONB in Hilbert space  $E$.
The theorems 4.9 and 5.1 imply the following statement. 

{\bf Theorem 6.1}. {\it Let ${\cal E},\, {\cal F}\in {\cal S}$. Then $\lambda _{\cal E}|_{{\cal R}_{\cal E}\bigcap {\cal R}_{\cal F}}=\lambda _{\cal F}|_{{\cal R}_{\cal E}\bigcap {\cal R}_{\cal F}}$. In particular, if the condition (\ref{crit}) is satisfied then ${\cal R}_{\cal E}={\cal R}_{\cal F}$ and $\lambda _{\cal E}=\lambda _{\cal F}$. If the condition (\ref{crit}) is violated then $\lambda _{\cal F}(A)=\lambda _{\cal E}(A)=0$ for every set $A\in {\cal R}_{\cal E}\bigcap {\cal R}_{\cal F}$.}

Consider the family  $\{ R_{\mathcal E},{\mathcal E}\in {\mathcal S}\}$ of rings of subsets of the space  $E$. Let ${\mathcal M}=\bigcup\limits_{{\cal E}\in {\cal S}}{\cal R}_{\cal E}$.
Let us define the function of a set $\lambda :\ {\mathcal M}\, \to \, [0,+\infty ]$ by the equality $\lambda (A)=\lambda _{\mathcal E}(A)\ \forall \ A\in R_{\mathcal E}$. The function  $\lambda $ is correctly defined since  $\lambda _{\mathcal E}(A)=\lambda _{\mathcal F}(A)$ for every $A\in R_{\mathcal E}\bigcap R_{\mathcal F}$ according to the theorem 6.1. Let $r$ be the ring generated by the family of sets  $\mathcal M$. Let's study the problem of continuation of the function  $\lambda$ from the collection of sets $\mathcal M$ onto the ring $r$ (see \cite{Tar}). 

Let's introduce the relation $\sim$ on the set $\cal S$ of ONB by the following way. ONB $\cal E$ and $\cal F$ are in the relation $\sim $ if the condition (\ref{crit}) is satisfied for ONB ${\cal E}$ and ${\cal F}$.

The relation $\sim$ in definition 6.2 is obviously reflexive.
According to the theorem 3.6 the relation $\sim$ is symmetric. Now we prove that the relation $\sim$ is transitive. Let's assume that the pairs of ONB $\cal E$, $\cal F$ and $\cal F$, $\cal G$ satisfy the relation $\sim$  (\ref{crit}). According to the theorem 6.1 the equalities ${\cal R}_{\cal E}={\cal R}_{\cal F}$ and ${\cal R}_{\cal F}={\cal R}_{\cal G}$ hold. Therefore, ${\cal R}_{\cal E}={\cal R}_{\cal G}$. Hence ONB  $\cal E$,  $\cal G$ satisfy the condition (\ref{crit}) according to the theorem 6.1. Hence the relation $\sim$ is transitive.

{\bf Definition 6.2}.
ONB $\cal E$ and $\cal F$ are called equivalent if they satisfy the condition (\ref{crit}). 

Let $\Sigma $ be a set of equivalence classes with respect to relation $\sim$: $\Sigma ={\cal S}/ \sim$. %
For every ONB ${\cal E}\in {\cal S}$ the space 
${\cal H}_{\cal E}=L_2(E,{\cal R}_{\cal E},\lambda _{\cal E},{\mathbb C})$ of square integrable with respect to the measure $\lambda _{\cal E}$ complex valued functions is introduced by the standard way (see  \cite{s16}). If $\{ {\cal E}\}\in \Sigma$ and ${\cal E'},{\cal E}''\in \{ {\cal E}\}$ then ${\cal H}_{{\cal E}'}={\cal H}_{{\cal E}''}$ according to the theorem 6.1 and definition of the spaces ${\cal H}_{\cal E},\, {\cal E}\in {\cal S}$. The symbol ${\cal H}_{\{ {\cal E}\} }$ notes the space ${\cal H}_{\cal E}$ for arbitrary choice of an ONB ${\cal E}\in \{ {\cal E}\}$.

Now we describe the ring generated by the family of subsets  ${\cal R}_{\cal E}\bigcup {\cal R}_{\cal F}$ for a pair of ONB ${\cal E},\, {\cal F}$ belonging to different classes $\{ {\cal E}\},\{{\cal F}\}\in \Sigma $. Moreover the sum of the spaces ${\cal H}_{\cal E}$ and ${\cal H}_{\cal F}$ will be defined.

The intersection of the rings ${\mathcal R}_{\mathcal E}$ and ${\mathcal R}_{\mathcal F}$ is the ring which is noted by the symbol ${\mathcal R}_{\mathcal E}\bigcap {\mathcal F}$. Then $\lambda _{\mathcal E}(A)=\lambda _{\mathcal F}(A)\equiv \lambda _{{\mathcal E} {\mathcal F}}(A)\;\;  \forall \ A\in  {\mathcal R}_{{\mathcal E}}\bigcap {\mathcal R}_{{\mathcal F}}$ according to the theorem 6.1. %
Then the space ${\mathcal H}_{{\mathcal E}\bigcap {\mathcal F}}=L_2(E, {\mathcal R}_{\mathcal E}\bigcap {\cal R}_{\mathcal F},\lambda _{{\mathcal E} {\mathcal F}}, \mathbb{C})$ is the subspace of Hilbert spaces ${\cal H}_{\cal E}$ and ${\cal H}_{\cal F}$.
Let ${\mathcal H}_{{\mathcal E}}^{\bot {\mathcal F}}$ and ${\mathcal H}_{{\mathcal F}}^{\bot {\mathcal E}}$ be orthogonal complements of the space ${\mathcal H}_{{\mathcal E}\bigcap {\mathcal F}}$ up to the spaces ${\mathcal H}_{\mathcal E}$ and ${\mathcal H}_{\mathcal F}$ respectively.
Then Hilbert space ${\mathcal H}_{{\mathcal E}{\mathcal F}}$ is defined as the direct sum of three orthogonal subspaces
\begin{equation}\label{summa}
{\mathcal H}_{{\mathcal E}{\mathcal F}}={\mathcal H}_{\mathcal E}^{\bot {\mathcal F}}\oplus {\mathcal H}_{{\mathcal E}\bigcap {\mathcal F}}\oplus {\mathcal H}_{\mathcal F}^{\bot {\mathcal E}}.
\end{equation}

{\bf Lemma 6.3}. {\it Let $\{ {\mathcal E}\},\, \{{\mathcal F}\}$ be the different equivalence classes of ONB in the space $E$. Let ${\cal E}\in \{ {\cal E}\},\, {\cal F}\in \{ {\cal F}\}$.
Then ${\cal H}_{\cal EF}={\cal H}_{\cal E}\oplus {\cal H}_{\cal F}$. Moreover, there is the shift-invariant measure   $\lambda _{{\cal E}{\cal F}}:\ {\cal R}_{{\cal EF}}\to [0,+\infty )$ such that ${\cal H}_{\cal E}\oplus {\cal H}_{\cal F}=L_2(E,{\cal R}_{{\cal EF}},\lambda _{{\cal E}{\cal F}},{\mathbb C})$. Here the ring  ${\cal R}_{{\cal EF}}$ is generated by the family of sets ${\cal R}_{\cal E}\bigcup {\cal R}_{\cal F}$}.

 If equivalence classes 
$\{ {\cal E}\},\{{\cal F}\}\in \Sigma $ are different then the equality  $\lambda _{\cal EF}(A)=0$ holds for any  $A\in {\mathcal R}_{{\mathcal E}} \bigcap {\cal R}_{{\mathcal F}}$ according to the theorem 6.1.
Hence the space ${\mathcal H}_{{\mathcal E}\bigcap {\mathcal F}}$ is trivial and ${\mathcal H}_{{\mathcal E}{\mathcal F}}={\mathcal H}_{\mathcal E}\oplus  {\mathcal H}_{\mathcal F}$. $\Box$

Thus every pair of Hilbert spaces ${\mathcal H}_{\{ {\mathcal E}\}}$ and ${\mathcal H}_{\{ {\mathcal F}\} }$ 
defines the Hilbert space  
\begin{equation}\label{raz}
{\mathcal H}_{\{ {\mathcal E}\}\{{\mathcal F}\}}= {\mathcal H}_{\{ {\mathcal E}\} }\oplus {\mathcal H}_{\{ {\mathcal F}\}}.
\end{equation}

Let ${\mathcal R}_{\{ {\mathcal E}\} \{ {\mathcal F}\} }$ -be the ring which is generated by the collection of sets ${\mathcal R}_{\{ {\mathcal E}\} }\bigcup {\mathcal R}_{\{ {\mathcal F}\} }.$  
The equality (\ref{raz}) defines the scalar product in the space ${\mathcal H}_{\{ {\mathcal E}\}\{{\mathcal F}\}}$. This scalar product defines (see  \cite{s16}) the continuation of the measures  $\lambda _{\{ \cal E\} },\, \lambda _{\{ \cal F\} }$ onto the shift-invariant measure $\lambda _{\{ \cal E\}\{ \cal F\} }:\ {\mathcal R}_{\{ {\mathcal E}\}\{  {\mathcal F}\} }\to [0,+\infty )$ according to the following condition.  The value of the measure $\lambda _{\{ \cal E\}\{ \cal F\} }$ on a set  $A\bigcap B$ is given by the equality
\begin{equation}\label{rin}
\lambda _{\{{\mathcal E} \}\{{\mathcal F}\} }(A\bigcap B)=(\chi _{A},\chi _B)_{{\cal H}_{{\cal EF}}} 
\end{equation}
for every sets  $A\in {\mathcal R}_{\{ \mathcal E\} },\, B\in {\mathcal R}_{\{ \mathcal F\} }$.
The value of the measure $\lambda _{\{ {\mathcal E} \} \{{\mathcal F}\}}$ on the other sets of the ring ${\mathcal R}_{\{{\mathcal E}\}\{ {\mathcal F}\}}$ is defined by the additivity condition. Therefore the function $\lambda _{\{{\mathcal E}\} \{{\mathcal F}\}}:\ {\mathcal R}_{\{{\mathcal E}\} \{{\mathcal F}\}}\to R$ is the finitely-additive measure. This measure is shift-invariant by the construction.
Moreover, if $A\in {\mathcal R }_{\{ {\mathcal E}\}},\, B\in {\mathcal R}_{\{ {\mathcal F}\} }$, then $\lambda _{\{{\mathcal E}\} \{{\mathcal F}\}}(A\bigcap B)=0$ and $\lambda _{\{{\mathcal E}\} \{{\mathcal F}\}}(A\bigcup B)=\lambda _{\{ {\cal E}\} }(A)+\lambda _{\{ {\cal F}\} }(B)$ according to  (\ref{raz}) and (\ref{rin}) .
$\Box$

Let's endow the linear hull ${\cal L}({\cal H}_{\{ {\cal E}\}},\, \{ {\cal E}\}\in \Sigma )$ with the euclidean norm of direct sum of Hilbert spaces. Let $u=\sum\limits_{k=1}^mv_k,$ where  $v_k\in {\cal H}_{\{ {\cal E}\}_k},\, k=1,...,m$. Then the intersection  ${\cal H}_{\{ {\cal E}\}_k},\, {\cal H}_{\{ {\cal E}\}_j},\, k\neq j,$ is trivial subspace according to lemma 6.3.
Hence the representation of an element $u\in {\cal L}({\cal H}_{\{ {\cal E}\}},\, \{ {\cal E}\}\in \Sigma ) $ in the form $u=\sum\limits_{k=1}^mv_k,$ где $v_k\in {\cal H}_{\{ {\cal E}\}_k},\, k=1,...,m$ is unique.  Therefore, the function $\| \cdot \|:\ {\cal L}({\cal H}_{\{ {\cal E}\}},\, \{ {\cal E}\}\in \Sigma )\to [0,+\infty )$ such that  
\begin{equation}\label{euclid}
\| u\|=(\sum\limits_{k=1}^m\| v_k\|_{{\{ {\cal E}\}}_k})^{1\over 2}, 
\end{equation}
for any vector $u\in {\cal L}({\cal H}_{\{ {\cal E}\}},\, \{ {\cal E}\}\in \Sigma ) $ of the form  $u=\sum\limits_{k=1}^mv_k,$  $v_k\in {\cal H}_{{\{ \cal E\}}_k},\, k=1,...,m$,  is the euclidean norm on the space  ${\cal L}({\cal H}_{\{ {\cal E}\}},\, \{ {\cal E}\}\in \Sigma )$. By the construction the linear hull ${\cal L}({\cal H}_{\{ {\cal E}\}},\, \{ {\cal E}\}\in \Sigma )$ with the norm (\ref{euclid}) are invariant both with respect to a shift of arguments of functions  $u\in {\cal L}({\cal H}_{\{ {\cal E}\}},\, \{ {\cal E}\}\in \Sigma )$ on any vector of the space  $E$ and with respect to a transformation of argument of functions $u\in {\cal L}({\cal H}_{\{ {\cal E}\}},\, \{ {\cal E}\}\in \Sigma )$ by any orthogonal transformation of the space $E$.
The completion of euclidean space ${\cal L}({\cal H}_{\{ { \cal E}\}},\, \{ {\cal E}\}\in \Sigma )$ by the norm  (\ref{euclid}) is the Hilbert space
$\mathbb H$.

Let  $\cal R$ be the minimal ring of subsets of the space $E$ containing any ring  ${\cal R}_{\cal E},\, {\cal E}\in {\cal S}$.

{\bf Theorem 6.4. } {\it There is the complete $\sigma$-finite locally finite finitely-additive  measure  $\lambda:\ {\cal R}\to [0,+\infty )$  which is invariant both with respect to shift on any vector of the space $E$ and with respect to any orthogonal transformation of the space $E$. The measure $\lambda $ is connected with the space $\mathbb H$ by the equality ${\mathbb H}=L_2(E,{\cal R},\lambda ,{\mathbb C})$. The measure $\lambda $ is not countably additive and $\lambda |{{\cal R}_{\{ {\cal E}\}}}=\lambda _{\{ {\cal E}\}}$ for every $\{ {\cal E}\}\in \Sigma$.}

{\bf Proof}.
Let $r$ be the ring of subsets of the space $E$ which is generated by the system of sets $\{ {\cal R}_{\{ { \cal E}\}},\, \{ {\cal E}\}\in \Sigma \}$. Hence the ring $r$ is generated by the semiring  
\begin{equation}\label{s-ring}
s=\{ A_0\backslash \bigcup\limits_{j=1}^NA_j,\ N\in {\mathbb N},\, A_0\in  {\cal R}_{\{ { \cal E}\}_0},\, A_j\in  {\cal R}_{\{ { \cal E}\}_j}\}.
\end{equation} 
Since the system of sets ${\cal R}_{\{ { \cal E}\}_0}$ is the ring then we can assume that $\{ {\cal E}\}_j\neq \{ {\cal E}\}_0$ for every  $j=1,...,N.$ Hence $ \lambda _{\cal E}(A_0\bigcap (\bigcup\limits_{j=1}^NA_j))=0$ according to the theorem 6.1. Thus we should define $\lambda (A_0\backslash \bigcup\limits_{j=1}^NA_j)=\lambda _{\cal E}(A_0)$ for any set  
$A_0\backslash \bigcup\limits_{j=1}^NA_j,\ j\in {\mathbb N},\, A_0\in  {\cal R}_{\{ { \cal E}\}_0},\, A_j\in  {\cal R}_{\{ { \cal E}\}_j}$ from the semiring (\ref{s-ring}). Then the function $\lambda :\ s\to [0,+\infty )$ is additive on the semiring (\ref{s-ring}). Moreover, this additive function satisfies the condition
$\lambda (A)=\| \chi _A\|_{\mathbb H}^2\ \forall \ A\in s$. Additive function $\lambda :\ s\to [0,+\infty )$ on the semi-ring $s$ admits the only one  
additive extension onto the additive function of a set $\lambda :\ r\to [0,+\infty )$ on the ring $r$.
Moreover, the measure $\lambda$ satisfy the condition $\lambda (A)=\| \chi _A\|_{\mathbb H}^2\ \forall \ A\in r$.

The semiring $s$ and the generated by this semiring ring $r$ are invariant with respect to both a shift on a vector of the space $E$ and an orthogonal mapping of the space $E$. The measure  $\lambda :\ r\to [0,+\infty )$ is both rotation and shift-invariant measure on the space $E$ too by its construction.

The function of a set $\lambda :\ r\to [0,+\infty )$ is finitely-additive by its construction. According to \cite{s16} (see also \cite{BusS}) the measure $\lambda $ takes zero values on a ball of the space $E$ with sufficiently small radius $\rho \in (0,{1\over {\sqrt 2}})$ (\cite{Busovikov}). Therefore the measure $\lambda $ is locally finite. Its $\sigma$-finiteness is the consequence of its locally finiteness and the separability of the space $E$. Moreover, since $\lambda (B_{\rho }(a))=0$ for any ball $B_{\rho }(a)=\{ x\in E:\ \|x-a\|_E<\rho \}$ where $a\in E,\ \rho={1\over 4}$ then the measure $\lambda$ is not countably additive according to the separability of the space $E$.

The measure $\lambda $ is not complete. Its completion $\lambda :\ {\cal R}\to [0,+\infty )$ is defined by the standard scheme by means of external and internal measures (see \cite{s16}).

According to the construction of measure $\lambda$  the following equality $\lambda |_{{\cal R}_{\{ {\cal E}\}}}=\lambda _{\{ {\cal E}\}}$ holds for every  $\{ {\cal E}\}\in \Sigma $.
Therefore, $L_2(E,{\cal R},\bar \lambda ,{\mathcal C}) \supset {\mathcal L}({\cal H}_{\{ { \cal E}\}},\, \{ {\cal E}\}\in \Sigma )$. 
Linear manifold ${\mathcal L}({\cal H}_{\{ { \cal E}\}},\, \{ {\cal E}\}\in \Sigma )$ is dense in the space $\mathbb H$ since the space ${\mathbb H}$ is defined as the completion of the euclidean space $({\mathcal L}({\cal H}_{\{ { \cal E}\}},\, \{ {\cal E}\}\in \Sigma ), \| \cdot \|)$.

The linear hull  ${\rm span}(\chi _A,\, A\in {\cal R})$ of the family of indicator functions of a sets from the ring $\cal R$ is the dense  linear manifold in the space $L_2(E,{\cal R},\lambda ,{\mathbb C})$ according to definition of this space.
Since the ring $\cal R$ is the completion of the ring $r$ with respect to measure $\lambda :\ r\to [0,+\infty )$ then the linear hull of the family of indicator functions ${\rm span}(\chi _A,\, A\in r)$ is dense in the space $L_2(E,{\cal R},\lambda ,{\mathbb C})$.

The ring $r$ is generated by the families of sets
$\{ {\cal R}_{\{ { \cal E}\}},\, \{ {\cal E}\}\in \Sigma \}$. Therefore, the linear manifold  ${\rm span}(\chi _A,\, A\in r_{\{ { \cal E}\}},\, \{ {\cal E}\}\in \Sigma )$ is dense in linear space  ${\rm span}(\chi _A,\, A\in r)$ equipped with the euclidean norm of the space $L_2(E,{\cal R}, \lambda ,{\mathbb C})$. 

Since the linear space ${\rm span}(\chi _A,\, A\in r_{\cal E})$ is the dense linear manifold in the space ${\cal H}_{\cal E}$ for every $\{ {\cal E}\}\in \Sigma $ then the linear manifold 
 $${\rm span}({\rm span}(\chi _A,\, A\in r_{\cal E}),\, \{ {\cal E}\}\in \Sigma )={\rm span}(\chi _A,\, A\in r_{\{ {\cal E}\}},\, \{ {\cal E}\}\in \Sigma )$$
dense in the linear space $${\mathcal L}({\cal H}_{\{ { \cal E}\}},\, \{ {\cal E}\}\in \Sigma )={\rm span}({\cal H}_{\{ { \cal E}\}},\, \{ {\cal E}\}\in \Sigma )$$ equipped with the norm $\| \cdot \|_{\cal H}$.

Therefore,
$L_2(E,{\cal R},\lambda ,{\mathbb C})=\mathbb H$ since the norm $\| \cdot \|_{\cal H}$ and the norm of the space $L_2(E,{\cal R},\lambda ,{\mathbb C})$ are coincide on the vectors of linear manifold ${\rm span}(\chi _A,\, A\in r_{\{ {\cal E}\}},\, \{ {\cal E}\}\in \Sigma )$ and this linear manifold is dense both in the space 
$L_2(E,{\cal R},\lambda ,{\mathbb C})$ and in the space  $\mathbb H$. $\Box$

The theorem 6.4 gives the orthogonal decomposition of the space 
${\mathbb H}$.
According to the lemma 6.3 the condition $\{ {\cal E}\} \neq \{ {\cal F}\}$ implies  ${\cal H}_{\{ {\cal E}\} }\bot {\cal H}_{\{ {\cal F}\} }$. Thus we obtain the following statement.

{\bf Corollary 6.5}.
$$
{\mathbb  H}=\oplus_{ \{ {\cal E}\} \in \Sigma}{\cal H}_{\{ {\cal E}\} }.
$$

In fact, the Hilbert space ${\mathbb H}$ is defined as the completion of linear hull of the family of spaces ${\mathcal H}_{\mathcal F},\, {\mathcal F}\in {\mathcal S}$ equipped with the scalar product (\ref{euclid}).

{\bf Corollary 6.6}. {\it Let ${\cal E},\, {\cal F}\in {\cal S}$. If bases  ${\cal E},\, {\cal F}$ satisfy the condition (\ref{crit}) then  ${\cal H}_{\cal E}={\cal H}_{\cal F}$. If the condition (\ref{crit}) is violated for bases ${\cal E},\, {\cal F}$ then ${\cal H}_{\cal E}\bot {\cal H}_{\cal F}$}.

{\bf Theorem 6.7}. 
{\it The shifts and rotation-invariant measure $\lambda :\ {\cal R}\to [0,+\infty )$ on the space $E$ from the theorem 6.4 admits the decomposition 
\begin{equation}\label{razlo}
\lambda =\sum\limits_{\{{\cal F} \}\in \Sigma}\nu _{\{ {\cal F}\}}
\end{equation}
into the sum of mutually singular shift-invariant measures $\nu _{\{ {\cal F}\}},\, \{{\cal F} \}\in \Sigma$. Here for every $\{ {\cal F}\}\in \Sigma $ the measure $\nu _{\{ {\cal F}\}}$ is given by the equality}
\begin{equation}\label{parti}
\nu _{\{ {\cal F}\}}(A)=\sup\limits_{B\in {\cal R}_{\cal F},\, B\subset A}\lambda _{\{ {\cal F}\}}(B),\, A\in {\cal R}.
\end{equation}

{\bf Proof}. In the proof of theorem 6.4 we show that if  $A\in \cal R$ then for any $\epsilon >0$ there are  $A_1'\in {\cal R}_{\{ {\cal E}_1\} },...,A_m'\in {\cal R}_{\{ {\cal E}_m\} }$ and $A_1''\in {\cal R}_{\{ {\cal E}_1\} },...,A_m''\in {\cal R}_{\{ {\cal E}_m\} }$ such that $\bigcup\limits_{j=1}^mA_j'\subset A\subset \bigcup\limits_{j=1}^mA_j''$ and $\lambda ((\bigcup\limits_{j=1}^mA_j'')\backslash (\bigcup\limits_{j=1}^mA_j'))<\epsilon $. This fact and the theorem 6.1 together imply the equality (\ref{razlo}) where measures $\nu _{\cal F}$ are defined by the (\ref{parti}). The mutually singularity of measures $\nu _{\{ {\cal E}\}},\, \nu _{\{ {\cal F}\}}$ under the assumption $\{ \cal E\}\neq \{ \cal F\}$ is the consequence of the decomposition 
(\ref{razlo}) and the theorem 6.1. $\Box$

{\bf Remark 6.8.}
If  ${\bf U}$ is the unitary operator in the space 
$E$ then $\nu _{\{ {\cal F}\}}({\bf U}A)=\nu _{\{ {\bf U}^{-1}{\cal F}\}}(A)$ for any set $A\in {\cal R}$ and for any class of ONB $\{ {\cal F}\} \in \Sigma$.

\section{Linear operators in the space $\mathbb H$ generated by the orthogonal transformation of argument} 

Let ${\bf U}(t),\,  t\in R$, be the one-parametric family of operators in the space 
${\mathbb H}$ which is given by the following way. Consider an ONB ${\cal E}=\{ e_j,\, j\in {\mathbb N}\}$ in the space $E$. Let $E_k={\rm span}(e_{2k-1},e_{2k}),\, k\in \mathbb N$, be a sequence of two-dimensional orthogonal subspaces of the space $E$. Let's consider the group  ${\bf \Lambda} (t),\, t\in {\mathbb R}$ of orthogonal transformations of the space $E$ such that, at first, the subspaces $E_k,\, k\in \mathbb N$ are invariant with respect to operators of this group and, at second, for every $k\in \mathbb N$ the restriction ${\bf \Lambda }(t)|_{E_k},\, t\in {\mathbb R}$ has the matrix 
$$
\left(
\begin{array}{llllll}
\cos (a_kt)&\sin (a_kt)\\
-\sin (a_kt)&\cos (a_kt)\\
\end{array}
\right).
$$
Here $\{ a_k\}$ is a sequence of real numbers.

Let   ${\bf U}(t),\,  t\in \mathbb R$, be a one-parametric family of operators in the space $E$ such that for every $t\in \mathbb R$ the operator  ${\bf U}(t)$ is given by the equality  
\begin{equation}\label{pri}
{\bf U}(t)u(x)=u({\bf \Lambda }(t)x),\, x\in E,\ u\in {\mathbb H}.
\end{equation}

Since the measure $\lambda :\ {\cal R}\to [0,+\infty )$ is invariant with respect to any orthogonal transformation of the space $E$ then the equality (\ref{pri}) defines the unitary operator  ${\bf U}(t)$ in the space $\mathbb H$ for every  $t\in \mathbb R$. One-parametric family of operators ${\bf U}(t), \, t\in \mathbb R$ forms the one-parametric unitary group in the space  
$\mathbb H$.

It is easy to check that the matrix of orthogonal operator ${\bf \Lambda }(t),\, t\neq 0,$ in an ONB $\cal E$ satisfies the condition (\ref{crit}) if and only if $\{ a_k\} \in \ell _1$.

{\bf Lemma 7.1}. {\it If $\{ a_k\}\notin \ell _1$ then one-parametric group ${\bf U}(t),\, t\in \mathbb R$, of unitary operators in the space  $\mathbb H$ is not strong continuous.}

{\bf Proof}.
Let $A\in {\cal R}_{\cal F}$ be a set such that $\lambda (A)=\lambda _{\cal F}(A)>0.$
Let's assume the opposite,  that the one-parametric group ${\bf U}(t),\, t\in \mathbb R$ is strong continuous. Then the function $({\bf U}(t)\chi _A,\chi _A)_{\mathbb H},\, t\in \mathbb R$, is continuous. 

But this function has the discontinuity point $t_0=0$ since $({\bf U}(t)\chi _A,\chi _A)_{\mathbb H}|_{t=0}=\lambda (A)>0,$ and $({\bf U}(t)\chi _A,\chi _A)_{\mathbb H}=0\, \forall \, t\neq 0$. In fact, $\chi _A\in {\cal H}_{\{ \cal E\}},\ {\bf U}(t)\chi _A\in {\cal H}_{\{ {\bf \Lambda (t)}({\cal E})\}}$. Since the orthogonal mapping ${\bf \Lambda }(t)$ is not satisfies the condition (\ref{crit}) then the subspaces ${\cal H}_{\cal E}$ and ${\cal H}_{\{ {\bf \Lambda (t)}({\cal E})\}}$ are orthogonal. 
The obtained contradiction proves the statement. $\Box$

{\bf Lemma 7.2}. {\it If $\{ a_k\}\in \ell _1$ then a subspace ${\cal H}_{\cal E}$ is invariant with respect to operators of one-parametric groups  ${\bf U}(t),\, t\in \mathbb R$ and the group  ${\bf U}(t)|_{{\cal H}_{\cal E}},\, t\in \mathbb R$ is strong continuous unitary group in the space  ${\mathcal H}_{\cal E}$.}

{\bf Proof}. The condition  $\{ a_k\}\in \ell _1$ implies that for every $t\in \mathbb R$ the matrix of orthogonal mapping  ${\bf \Lambda} (t)$ in the basis $\cal E$ satisfies the condition (\ref{crit}). Therefore, ONB ${\cal E}$ and ${\bf \Lambda }(t)({\cal E})$ are equivalent and hence  ${\bf U}(t){\cal H}_{\cal E}={\cal H}_{\cal E}\ \forall \ t\in \mathbb R$. 

Let's fix a number $\epsilon >0$. Let $\phi \in {\cal H}_{\cal E}$.
For every $m\in \mathbb N$ the equality $\{a_k\}=\{ a_k'(m)\}+\{a_k''(m)\}$ holds. Here $\{ a_k'(m)\}=\{ a_1,...,a_m,0,...\},\ \{a_k''(m)\}=\{ 0,...,0,a_{m+1},...\}$. Hence ${\bf \Lambda}(t)={\bf \Lambda }_m''(t)\circ {\bf \Lambda }_m'(t),\, t\in {\mathbb R}$ and ${\bf U}(t)|_{{\cal H}_{\cal E}}={\bf U }''_m(t)|_{{\cal H}_{\cal E}}\circ {\bf U }'_m(t)|_{{\cal H}_{\cal E}},\, t\in {\mathbb R}$. Then for any $t>0$ there is a number $m\in \mathbb N$ such that $\sup\limits_{\tau \in [0,t]}\| {\bf U}''(\tau )_m\phi -\phi \|_{{\cal H}_{\cal E}}<\epsilon $.
In fact, for every $\Pi \in {\cal K}_{\cal E}$ we have $\| {\bf U}''(t)_m\chi _{\Pi }-\chi _{\Pi }\|_{\mathbb H}^2=2\lambda (({\bf \Lambda }''_m(t)\Pi )\backslash \Pi)$ and $\lim\limits_{m\to \infty}\sup\limits_{\tau \in [0,t]}\lambda (({\bf \Lambda }''_m(\tau )\Pi )\backslash \Pi)=0$ according to lemma 4.7. The strong continuity of the group ${\bf U}'_m(t)|_{{\cal H}_{\cal E}},\, t\in \mathbb R$ is the consequence of the decomposition  ${\cal H}_{\cal E}=L_2({\mathbb R}^{2m})\oplus {\cal H}_{{\cal E}''}$ (see \cite{BusS}) and the strong continuity of the group of orthogonal transformations of the argument of a square integrable function on a finite-dimensional euclidean space ${\mathbb R}^{2m}$. $\Box$

{\bf Remark 7.3}. Let $\{ a_k\}\in \ell _1$. Then a subspace ${\cal H}_{\cal F}$ can be not invariant with respect to operators of the one-parametric group ${\bf U}(t),\, t\in \mathbb R$ for every  ${\cal F}\in {\cal S}$.  Moreover, the group ${\bf U}(t)|_{{\cal H}_{\cal E}},\, t\in \mathbb R$, of unitary operators in the space ${\mathcal H}$ can be discontinuous. This fact is shown by the following example.

Let the operator ${\bf \Lambda} (t)$ in the ONB $\cal E$ has the matrix
$$
\| {\bf \Lambda }(t)\|_{\cal E}=\left(
\begin{array}{llllll}
\cos (at)&\sin (at)&0&0&0&\ldots\\
-\sin (at)&\cos (at)&0&0&0&\ldots\\
0&0&1&0&0&\ldots\\
0&0&0&1&0&\ldots\\
0&0&0&0&1&\ldots\\
\ldots&\ldots&\ldots&\ldots&\ldots&\ldots
\end{array}
\right)
$$
The operator ${\bf I}-{\bf \Lambda}(t)$ is trace class operator and the statement of the lemma 7.2 takes place. Let's consider the ONB $\cal F$ which is the image of the ONB $\cal E$ under the action of the orthogonal mapping $\bf V$ of the space $E$ with the following matrix in ONB $\cal E$ 
$$
\| {\bf V}\|_{{\cal E}}=\left(
\begin{array}{lllll}
1&0&0&0&\ldots\\
0&c_{2,2}&c_{2,3}&c_{2,4}&\ldots\\
0&c_{3,2}&c_{3,3}&c_{3,4}&\ldots\\
0&c_{4,2}&c_{4,3}&c_{4,4}&\ldots\\
\ldots&\ldots&\ldots&\ldots&\ldots
\end{array}
\right).
$$
Rows of the matrix $
\| {\bf V}\|_{{\cal E}}$ forms the orthonormal system in the space $\ell _2$. There is the choice of this orthonormal system such that the following condition 
\begin{equation}\label{rows}
\sum\limits_{k=1}^{\infty }|c_{2,k}|=+\infty 
\end{equation} 
holds. Let's consider the  orthogonal mapping $\bf V$ such that the condition (\ref{rows}) is satisfied.
Then the matrix of the orthogonal operator ${\bf \Lambda }(t)$ in the ONB $\cal F$ is
$$
\| {\bf \Lambda}(t)\|_{\cal F}=\| {\bf V}\|_{{\cal E}}^{-1} \| {\bf \Lambda }(t)\|_{\cal E}\| {\bf V}\|_{{\cal E}}=\left(
\begin{array}{lllll}
\cos (at)&c_{2,2}\sin (at)&c_{2,3}\sin (at)&c_{2,4}\sin (at)&\ldots\\
-c_{2,2}\sin (at)&.&.&.&\ldots\\
-c_{2,3}\sin (at)&.&.&.&\ldots\\
-c_{2,4}\sin (at)&.&.&.&\ldots\\
\ldots&\ldots&\ldots&\ldots&\ldots\\
\end{array}
\right).
$$
Hence the condition (\ref{crit}) is not satisfied for the bases  $\cal F$ and ${\cal F}'={\bf \Lambda }(t)\cal F$. Thus, the bases $\cal E$ and ${\cal E}'={\bf \Lambda }(t){\cal E}$ are equivalent in the sence of the definition 6.2, but the bases ${\cal F}$ and ${\cal F}'={\bf \Lambda }(t){\cal F}$ are not equivalent.  Hence the subspaces  ${\cal H}_{\cal F}$ and ${\bf U }(t){\cal H}_{\cal F}={\cal H}_{{\cal F}'}$ are orthogonal in the space $\cal H$ for any $t\neq 0$ according to the corollary 6.6. Therefore, the group  ${\bf U}(t),\, t\in \mathbb R$, of unitary operators (\ref{pri}) in the space $\mathbb H$ is discontinuous since the function $t\to (u,{\bf U}(t)u)_{\mathbb H},\, t\in \mathbb R$, is discontinuous for any nontrivial vector $u\in {\cal H}_{\cal F}.$

{\bf Remark 7.4}. The condition (\ref{crit}) on the orthogonal mapping $\bf \Lambda$ in the space  $E$ and ONB $\cal E$ can not be considered as the condition on the operator $\bf \Lambda$ only. In particular, the condition  (\ref{crit}) is not the consequence of the belonging of the operator  ${\bf I}-{\bf \Lambda}$ to the space of trace class operators. This fact is shown by the example in the remark 7.3. In fact, the operator  ${\bf \Lambda }(t)-{\bf I}$ is trace class operator (as well as the operator ${\bf V}^{-1}{\bf \Lambda }(t){\bf V}-{\bf I}$). Nevertheless, the operator ${\bf \Lambda }(t)$ and the basis $\cal E$ satisfy the condition  (\ref{crit}), but the operator ${\bf \Lambda }(t)$ and the basis ${\cal F}={\bf V}\cal E$ are not satisfy the condition (\ref{crit}). 
This fact is not the contradiction for the rotation invariance of the measure $\lambda$. Because the bases  ${\cal F}={\bf V}{\cal E}$ and ${\cal F}'={\bf \Lambda }(t){\cal F} 
$ are not equivalent, but the bases  ${\cal F}
$ and  ${\cal G}
={\bf V}{\bf \Lambda }(t){\bf V}^{-1}{\cal F}$ are equivalent.

{\bf Remark 7.5}. The  unitary group of operators in the space  $\mathbb H$ generated by the group of orthogonal mappings (\ref{pri}) in the space $E$ has the strong continuity property describing by lemmas 7.1, 7.2 and remark 7.3. This properties are similar to the strong continuity  properties of the unitary group of operators in the space $\mathbb H$ generated by the group of shifts of  argument according to the formula  
\begin{equation}\label{shi}
{\bf S}_{th}u(x)=u(x-th),\, t\in \mathbb R.
\end{equation}
The group of unitary operators (\ref{shi}) has the the family of invariant subspaces ${\cal H}_{\cal F},\ {\cal F}\in {\cal S}$ and the restriction ${\bf S}_{th}|_{{\cal H}_{\cal F}},\, t\in \mathbb R$, is strong continuous group in the space ${\cal H}_{\cal F}$ if and only if  $\{ (h,f_k)\}\in l_1 $ (see \cite{SZ}).

\section{Conclusion}

In this paper we construct the finitely-additive measure $\lambda$ on an infinite-dimensional real Hilbert space such that this measure is shift  and rotation invariant. Moreover, the introduced measure is locally finite and $\sigma$-finite. But it is not countably additive and Borel measure. The decomposition of the measure $\lambda$ into the sum of mutually singular shift-invariant measures is obtained.

Finitely-additive measures have some peculiar properties with respect to countably additive ones. Nevertheless they are an effective mathematical technique to study the theory of quantum structures
\cite{Brat-Rob}, to describe the set of singular quantum states
\cite{S11, AS16} and in the theory of quantum measurements of observables with the continuous spectrum
\cite{Srinivas}.
Finitely-additive measures (in particular, Feynman measures) are a foundation to the representation of the solution of initial-boundary value problems for an evolutionary differential equations by means of functional integrals
\cite{Smolyanov, S-Tsoy}. Shift and rotation-invariant measures on a Hilbert space are natural for the studying a random walks and diffusions in this space since the above invariance presents the properties of homogeneity and isotropy of the space.
Shift-invariant measure $\lambda $ on a Hilbert space $E$ defines the connection between the Gaussian random vector in a Hilbert space $E$ and  the self adjoint operator $\bf \Delta $ in the space $L_2(E,{\cal R},\lambda ,{\bf C})$ (see \cite{SZ}). 
Both the above self-adjoint  operators $\Delta $ and Levi-Laplacians (\cite{Volkov}) are 
infinite-dimensional analogs of Laplace operator.

Finitely-additive measures have the peculiar properties such as the absence of statements of the majorizable convergence theorem and  the Fubini theorem. Nevertheless, a lot of the constructions of functional analysis hold for the functions of infinite-dimensional argument, which is integrable with respect to shift-invariant finitely-additive measure. Namely, the analogs of Sobolev spaces and the space of smooth functions are defined, the embedding theorems for the scales of Sobolev spaces  and the space of smooth functions are obtained \cite{BusS}, the theorems on the traces are proved \cite{BusS2}, the boundary value problems for the Poisson equation in an infinite-dimensional domains are posed and the variational approaches for their solution are obtained \cite{BusS3}.

Shift invariant finitely-additive measures and corresponding functional spaces are used for the description of diffusion processes in a Hilbert space and to describe the quantum dynamics for quantizations of infinite-dimensional Hamiltonian systems
(\cite{s16, BusS, OSZ}).

The invariance of a measure on euclidean space with respect to orthogonal mappings has applications in differential equations, statistical and quantum mechanics. For example, the compositions of independent random orthogonal mappings in the finite-dimensional euclidean space and corresponding diffusion equation on the finite-dimensional sphere are studied in  
\cite{SSZ} by means of invariant with respect to orthogonal mappings Lebesgue measure.

\bigskip

{\large {\bf Acknowledgments}}.
The author thank Professor I.V. Volovich for the interest to this work
and important discussion of the analysis of finitely-additive measures and its applications in the general problems of invariant measures on groups, in the theory of Hamiltonian systems and in the quantum mechanics.

\end{document}